\documentclass[
]
{amsart}
\usepackage{amsmath,amssymb,amsthm,eucal,amscd}

\newtheorem{thm}{Theorem}[section]
\newtheorem{cor}[thm]{Corollary}

\newtheorem{lem}[thm]{Lemma}
\newtheorem{nota}[thm]{Notation}
\newtheorem{prop}[thm]{Proposition}
\newtheorem{as}[thm]{Assumption}

\theoremstyle{definition}
\newtheorem{df}[thm]{Definition}

\newtheorem{rem}[thm]{Remark}

\numberwithin{equation}{section}

\begin{document}

\def\<{\langle}
\def\>{\rangle}
\def\({\<\<}
\def\){\>\>}
\def\tri{| \hskip-0.02in |\hskip-0.02in |}
\def\1{{{1\hskip-2.0pt \rule{3pt}{0.2pt}\hskip-2.0pt
\rule{0.7pt}{6.6pt}
\hskip-2.2pt{\raise6.8pt\hbox{\rule{3pt}{0.2pt}}}\hskip0.8pt}}}
\def\1{{{{\Large \text{1}}\hskip-2.5pt \rule{3pt}{0.2pt}\hskip-2.0pt
\rule{0.7pt}{7.5pt}
\hskip-2.2pt{\raise7.6pt\hbox{\rule{3pt}{0.2pt}}}\hskip0.8pt}}}
\def\lnorm {\left|}
\def\rnorm {\right|}
\def\ad{\mathop{\rm ad}\nolimits}

\def\a{\mathfrak{A}}
\def\A{\mathcal{A}}
\def\C{\mathcal C}
\def\D{\mathcal{D}}
\def\e{\varepsilon}
\def\f{\mathfrak{F}}
\def\F{\mathcal{F}}
\def\g{\mathfrak g}
\def\G{\g_{\infty}}
\def\GG{G_{\infty}}
\def\H{\mathcal H}
\def\J{\mathcal J}
\def\l{\lambda}
\def\M{\widetilde{M}}
\def\n{\nabla}
\def\o{\omega}
\def\O{\Omega}
\def\P{\mathfrak{P}}
\def\p{\Phi}
\def\r{\vec{r}}
\def\R{\mathbb{R}}
\def\s{\vec{s}}
\def\X{\mathbf{X}}
\def\y{\delta}
\def\z{\frac}

\def\so{\mathfrak{s}\mathfrak{o}_{HS}}
\def\gl{\mathfrak{g}\mathfrak{l}_{HS}}

\def\CM{G_{CM}}

\def\ll#1{\label{#1}%
}

\def\zoxit#1
{\leavevmode\vbox{\hbox to 0pt{\hss
\raise1.8ex
\vbox to 0pt
{\vss\hrule\hbox
{\vrule\kern.75pt\vbox{\kern.75pt
\hbox{\tiny #1}
\kern.75pt}\kern.75pt\vrule}\hrule}}}}
\def\zb#1#2{\begin{#1}\ll{#2} 
}

\title[Riemannian geometry of infinite-dimensional Lie groups]
{Hilbert-Schmidt groups as infinite-dimensional Lie groups and their Riemannian geometry}
\author{Maria Gordina}
\address{Department of Mathematics,
University of Connecticut, Storrs, CT 06269, U.S.A. }

\email{gordina@math.uconn.edu}
\date{\today}

\footnote{Research supported by the NSF Grant DMS-0306468.}

\begin{abstract}
We describe the exponential map from an infinite-dimensional Lie algebra to an infinite-dimensional group of operators on a Hilbert space. Notions of differential geometry are introduced for these groups. In particular, the Ricci curvature, which is understood as the limit of the Ricci curvature of finite-dimensional groups, is calculated. We show that for some of these groups the Ricci curvature is $-\infty$. 
\end{abstract}

\keywords{infinite-dimensional groups, Lie groups and Lie algebras, exponential map, Ricci curvature}

\maketitle

\renewcommand{\contentsname}{Table of Contents}

\tableofcontents

\section{Introduction}\ll{S:INTRO}

The results of this article are inspired by our previous study of heat kernel measures on infinite-dimensional Lie groups in \cite{Gor00-1}, \cite{Gor00-2}, \cite{Gor02}, \cite{Gor03-2}. The main tool in these papers was the theory of stochastic differential equations in infinite dimensions. The present paper, however, is entirely non probabilistic. 
It is organized as follows: in Sections \ref{S:compl} and \ref{S:CM} we discuss the exponential map for a certain class of infinite-dimensional groups, and in Sections \ref{S:geometry} and \ref{S:Ricci} we introduce notions of Riemannian geometry for Hilbert-Schmidt groups and compute the Ricci curvature for several examples.

\subsection{Motivation: Wiener measures and geometry} In our previous papers we were concerned with a pair of infinite dimensional Lie groups, $\CM \subset G_W$, related to each other in much the same way that the Cameron-Martin Hilbert space, $H_1([0,1])$, is related to Wiener space, $C_*([0,1])$: it is well understood that the geometry of the Hilbert space $H_1([0,1]) \subset  C_*([0,1])$ controls Wiener measure on $C_*([0,1])$, even though $H_1([0,1])$ is a subspace of Wiener measure zero. In the papers \cite{Gor00-1}, \cite{Gor00-2}, \cite{Gor02}, \cite{Gor03-2} we constructed an analog of Wiener measure on an infinite dimensional group $G_W$ as the ``heat kernel" (evaluated at the identity of $G_W$) associated to the Laplacian on the dense subgroup $\CM$. To this end one must choose an inner product on the Lie algebra, $\g_{CM}$, of $\CM$ in order to introduce a left invariant Riemannian metric on $\CM$. The Lie algebra $\g_{CM}$ determines a Laplacian, whose heat kernel measure actually lives on the larger group $G_W$. It can be shown that in general $\CM$ is a subgroup of measure zero. More features of this group have been discussed in \cite{Gor00-2}, \cite{Gor02}, \cite{Gor03-2}. In these papers we constructed the heat kernel measure by probabilistic techniques. We used a stochastic differential equation based on an infinite dimensional Brownian motion in the tangent space at the identity of $G_W$ whose covariance is determined by the inner product on $\g_{CM}$. 

Just as in the case of the classical Cameron-Martin space, where the Sobolev norm on $H_1([0,1])$ is much stronger than the supremum norm on $C([0,1])$, so also the norm on $\g_{CM}$ must be much stronger than the natural norm on the tangent space at the identity of $G_W$  in order for the heat kernel measure to live on $G_W$. If $G_W$ is simply the additive group, $C_*([0,1])$, and $\CM$ is the additive group $H_1([0,1])$, then this heat kernel construction reproduces the classical Wiener measure.

In this paper we address questions relating to the geometry of $\CM$, with a view toward eventual application to further understanding of the heat kernel measure on $G_W$. On a finite dimensional Riemannian manifold properties of the heat kernel measure are intimately related to the Ricci curvature of the manifold, and in particular to lower bounds on the Ricci curvature \cite{StroockBook}. We are going to compute the Ricci curvature for several classes of infinite dimensional groups, in particular, for those groups $\CM$ whose heat kernel measure on $G_W$ we have already proven the existence of in \cite{Gor00-1}, \cite{Gor00-2}, \cite{Gor02}, \cite{Gor03-2}. Our results show that the Ricci curvature is generally not bounded below, even in the cases when we were able to construct the heat kernel measure on $G_W$ (e.g. for the group $SO_{HS}$). One of the implications of our results is that the methods used to prove quasi-invariance of the heat kernel measure in the finite-dimensional case are not applicable for the settings described in our earlier papers. We also compute the Ricci curvature of groups $G_W$. All these groups are Hilbert-Schmidt groups which are described below.

\subsection{Hilbert-Schmidt groups as Lie groups and their Riemannian geometry } Denote by $HS$ the space of Hilbert-Schmidt operators on a real separable Hilbert space $H$. Let $B(H)$ be the space of bounded operators on $H$, and let $I$ be the identity operator. Denote by $GL(H)$ the group of invertible elements of $B(H)$. The general Hilbert-Schmidt group is $GL_{HS}=GL(H)\cap (HS+I)$. 
In the papers \cite{Gor00-1}, \cite{Gor00-2}, \cite{Gor02}, \cite{Gor03-2} we proved the existence and some basic properties of the heat kernel measures on certain classical subgroups of $GL_{HS}$, namely,  $SO_{HS}$ and $Sp_{HS}$. The Lie algebras of these groups are closed subspaces of $HS$ in the Hilbert-Schmidt norm. In the setting described above $GL_{HS}$, $SO_{HS}$ and $Sp_{HS}$ are examples of the group $G_W$. 

But the corresponding Cameron-Martin subgroups are, of necessity, only dense subgroups of $G_W$. They are determined by their tangent space, $\g_{CM}$, at the identity. In order to get the corresponding heat kernel measure to live on $GL_{HS}$, $SO_{HS}$ or $Sp_{HS}$, respectively, the tangent space $\g_{CM}$ must be given a Hilbert norm $|\cdot|$ which is much stronger than the Hilbert-Schmidt norm. The result is that the commutator bracket of operators may not be continuous in this norm. That is, $|AB-BA|\leqslant C|A||B|$ may fail for any constant $C$ as $A$ and $B$ run over $\g_{CM}$. Consequently $\g_{CM}$ may not really be a Lie algebra. Rather, the commutator bracket may be only densely defined as a bilinear map into $\g_{CM}$. In Section \ref{S:Lie} we will give a class of examples of groups contained in $GL_{HS}$ such that $\g_{CM}$ is not closed under the commutator bracket and in this sense is not a Lie algebra. But in most of the examples we consider this is not the case.

In this paper we are going first to address the problem of the relation of the tangent space, $T_I(\CM)$, to the Lie algebra structure of some dense subspace of  $T_I(\CM)$. As in the case of classical Wiener space, where polygonal paths play a central technical role in the work of Cameron and Martin, so also it seems to be unavoidable, for the purposes of \cite{Gor00-1}, \cite{Gor00-2}, \cite{Gor02}, \cite{Gor03-2}, to make use of a group $G \subset \CM$ which is in some sense dense in $\CM$ and which is itself a union of an increasing sequence of finite dimensional Lie groups: $G=\bigcup\limits_{n=1}^{\infty} G _n$, $G_n \subseteq G_{n+1}$. In Sections \ref{S:compl} and  \ref{S:CM} all of our groups will be taken to be subgroups of $GL(H)$ rather than $GL_{HS}$. Then the Lie algebra of $G_n$, $\g_n$, is a finite dimensional subspace of $B(H)$, and is closed under the commutator bracket. Let $\g=\bigcup\limits_{n=1}^{\infty} \g _n$. Then $\g$ is a also a Lie algebra under the commutator bracket. But $\g$ cannot be complete in any norm in the infinite dimensional case. 

Our goal in Sections \ref{S:compl} and  \ref{S:CM} is to study the completion, $\G$, of $\g$ in some (strong) Hilbert norm. In these sections we will assume that the completion $\G$ actually embeds into $B(H)$ and that the commutator bracket is continuous in this norm. We will characterize the group, $\CM$, generated by $\exp(\G)$ in this case, and show that the exponential map covers a neighborhood of the identity. These groups are examples of so called Baker-Campbell-Hausdorff Lie groups (e.g. \cite{Gl02-1}, \cite{Gl02-2}, \cite{Gl02-3}, \cite{R97}, \cite{MR95}, \cite{R96}, \cite{V72}, \cite{W98}). Let us mention here that the question of whether the exponential map is a local diffeomorphism into an infinite-dimensional Lie group has a long history. Our treatment is different in two major aspects. The first one is the choice of an inner product on $\g$ and corresponding norm on $\g$. As we mentioned earlier the heat kernel analysis on $GL_{HS}$ forces us to choose an inner product on $\g$ which is different from the Hilbert-Schmidt inner product. We will assume that the commutator bracket is continuous on $\G$ in the extended norm $|\cdot|$, namely, $|[x,y]|\leqslant C|x||y|$, where the constant $C$ is not necessarily $2$. In most results on Banach-Lie groups this constant is assumed to be $2$ (e.g. \cite{Har72}, \cite{LT65}). Quite often the underlying assumption is that $\g$ is a Banach algebra, and that $\g$ is complete. None of these is assumed in our case since we wish to deal with examples without these restrictions.  

In Sections \ref{S:geometry} and \ref{S:Ricci} we will compute the Ricci curvature in two major cases: when the norm on $\g$ is the Hilbert-Schmidt norm, and when it strongly dominates the Hilbert-Schmidt norm. 
Then we will examine how the lower bound  of the Ricci curvature depends on the choice of the strong Hilbert norm. We will extend Milnor's definitions of curvature on Lie groups \cite{Milnor76} to our infinite dimensional context for this purpose. Our results show that the Ricci curvature is generally not bounded below, and in some cases is identically minus infinity.

\subsection{Historical comments} We give references to the following mathematical literature addressing different features of infinite-dimensional Lie groups and exponential maps. Our list is certainly not complete, since the subject has been studied for many years. There are several reviews on the subject (e.g. \cite{Har72}, \cite{KM97}, \cite{Mi84}, \cite{Ne96}, \cite{R02}). Possible non-existence of an exponential map is addressed in \cite{Le92},  \cite{L95}, \cite{OMYK83}. Super Lie algebras of super Lie groups have been studied in \cite{Ci89}, \cite{Cz89}. Direct and inductive limits of finite-dimensional groups and their Lie group structures have been discussed in \cite{ARS85}, \cite{Gl03}, \cite{Gl02-3}, \cite{NRW01}, \cite{NRW91}, \cite{NRW94}, \cite{NRW93}, \cite{Neeb98}. Our main tool in proving that the exponential map is a local diffeomorphism is the Baker-Campbell-Dynkin-Hausdorff formula. Note that in the terminology of Banach-Lie groups it means that we prove that the groups we consider are Baker-Campbell-Hausdorff Lie groups (e.g. \cite{Gl02-1}, \cite{Gl02-2}, \cite{Gl02-3}, \cite{R97}, \cite{MR95}, \cite{R96}, \cite{V72}, \cite{W98}).  
Some of these articles also address the issue of completeness of the space over which a Lie group is modeled. We show that under Completeness Assumption \ref{A:compl} and Assumption \ref{E:norm} the Cameron-Martin group $\CM$ is complete in the metric induced by the inner product on $\g$. 
In addition, in Section \ref{S:Lie} we show that a natural completion of the infinite-dimensional Lie algebra $\g$ is not a Lie algebra in general without these assumptions.

One of the main contributors to the field of connections between differential geometry and stochastic analysis is P. Malliavin, who wrote a survey on the subject in \cite{M01}.
A book on stochastic analysis on manifolds has been written by E.Hsu \cite{HsuBook}. In conclusion we refer to works of B. Driver, S. Fang, D. Freed, E. Hsu, D. Stroock, I. Shigekawa, T. Wurzbacher et al dealing with infinite-dimensional Riemannian geometry and its applications to stochastic analysis (\cite{Dr92}, \cite{Dr94}, \cite{Driver95}, \cite{DriverLohr}, \cite{DriverHu}, \cite{Dr97}, \cite{Dr97-2}, \cite{Driver98}, \cite{Driver03}, \cite{Fang99}, \cite{Fang04}, \cite{Freed88}, \cite{HsuIAS}, \cite{Hsu99}, \cite{Hsu95}, \cite{Hsu93}, \cite{SegalBook}, \cite{Shig97}, \cite{SperaWurz00}, \cite{ShigTanig}, \cite{StroockBook}, \cite{Tan96}). They mostly concern loop groups, path spaces, their central extensions etc. In these cases the Riemannian geometry on the infinite-dimensional manifold is induced by the geometry of the space in which the loops or paths lie. The situation we consider is quite different.  

{\bf Acknowledgment.} I thank B.~Driver and P. Malliavin who asked me about the Ricci curvature of the Hilbert-Schmidt groups, and especially B.~Driver who suggested to use J. Milnor's \cite{Milnor76} results for finite-dimensional Lie groups. Section~\ref{S:CM} benefited greatly from my discussions with A.Teplyaev about the exponential map. Finally, I am very grateful to L.Gross for carefully reading the manuscript and suggesting significant improvements to the text.

\section{Completeness Assumption and the definition of the Cameron-Martin group}\ll{S:compl}

In this section we study a Lie group associated with an infinite-dimensional Lie algebra. The results of this section are not restricted to the Hilbert-Schmidt operators. We begin with an informal description of the setting. Let $\g$ be a Lie subalgebra of $B(H)$, the space of bounded linear  operators on a separable Hilbert space $H$. The group under consideration is a subgroup of $GL(H)$, the group of invertible elements of $B(H)$. The space $B(H)$ is the natural (infinite-dimensional) Lie algebra of $GL(H)$ with the operator commutator as the Lie bracket. 

We assume that $\g$ is equipped with a Hermitian inner product $(\cdot,\cdot)$, and the corresponding norm is denoted by $|\cdot |$. In an infinite-dimensional setting $\g$ might not be complete. We will always work with the situation when $\g$ has a completion which is a subspace of $B(H)$. But as we will see in Section \ref{S:Lie}, the most natural candidate for such a completion of $\g$ might not be closed under taking the Lie bracket. Similarly, when we look at the Lie group corresponding to the infinite-dimensional Lie algebra $\g$, it might not be complete in the metric induced by the inner product on the Lie algebra. 

The infinite-dimensional Lie algebra $\g$ is described by finite-dimensional approximations. Let $G_1 \subseteq G_2 \subseteq ... \subseteq G_n \subseteq ...\subseteq B(H)$ be a sequence of connected finite-dimensional Lie subgroups of $GL(H)$. 
Denote by $\g _n\subseteq B(H)$ their Lie algebras. We will consider the Lie algebra $\g=\bigcup\limits_{n=1}^{\infty} \g _n$. 

\begin{as}[Completeness Assumption]\ll{A:compl}
Throughout this section we assume  that there is a subspace $\G$ of $B(H)$ such that the Lie algebra $\g$ is contained in $\G$ and the given inner product $(\cdot,\cdot)$ on $\g$ extends to $\G$, which is complete with respect to this inner product. We will abuse notation by using $(\cdot,\cdot)$ to denote the extended inner product on $\G$ and by $|\cdot |$ the corresponding norm. We also assume that $\g$ is dense in $\G$ in the norm $|\cdot |$.
\end{as}

We will discuss this assumption in more detail in Section \ref{S:Lie} in the case of the Hilbert-Schmidt groups. In particular, we will show that $\G$ is not a Lie algebra in some cases of particular interest.

\begin{nota}Let $C^1_{CM}$ denote the space of paths $g:[0,1]\rightarrow GL(H)$ such that  

\begin{enumerate}

\item $g(s)$ is continuous in the operator norm,
\item
$\dot g=\frac {dg}{ds}$ exists in $B(H)$ equipped with the operator norm, 
\item $\dot g$ is piecewise continuous in the operator norm, 
\item $g'=g^{-1} \dot g$ is in $\G$, and $g'$ is piecewise continuous in the norm $|\cdot|$. 
\end{enumerate}
Let 

\[
d(y, z)=
\inf\limits_{g}\{\int\limits_0^1\lnorm g^{-1} \dot g \rnorm ds\},
\] 
where $g$ runs over $C^1_{CM}$ with $g(0)=y,g(1)=z$. We set $d(y, z)=\infty$ if there is no such path $g$. Note that $d$ depends on the norm $|\cdot |$ on $\g$. 
\end{nota}

\begin{nota}
$\CM=\{x \in B(H): d(x,I) < \infty \}$. 
\end{nota}

\begin{prop}\ll{St:group}
$\CM$ is a group,  and $d$ is a left-invariant metric on $\CM$.
\end{prop}

\begin{proof}
The proof for the first part is the same as for any finite-dimensional Lie group. In particular, if $f: [0,1]\rightarrow \CM, f(0)=x,f(1)=y$, $x, y \in \CM$, then  $h(s)=y^{-1}f(s)$ is a curve connecting $y^{-1}x$ and $I$, and 
$\lnorm h^{-1}\dot h \rnorm=
\lnorm f^{-1}\dot f \rnorm$. 
Therefore $d(y^{-1}x,I)=d(x,y)$, 
thus $y^{-1}x \in \CM$.
\end{proof}

\begin{df}\ll{D:CM}
$\CM$ is called \emph{the Cameron-Martin group}. 
\end{df}

\section{The Cameron-Martin group and the exponential map}\ll{S:CM}

\begin{as}[Continuity Assumption on the Lie bracket]\ll{E:norm}Throughout this section we assume that $\G$ is closed under taking the commutator bracket, and that the commutator bracket is continuous on $\G$, that is, there is $C>0$ such that

\[
\lnorm [x,y]\rnorm \leqslant C \lnorm x\rnorm \lnorm y\rnorm 
\]
for all $x,y \in \G$.
\end{as}

\begin{rem}
Assumption \ref{E:norm} is satisfied for any Banach algebra with $C=2$. In particular, it holds for the operator norm $\|\cdot \|$ and the Hilbert-Schmidt norm $|\cdot |_{HS}$, but these are not the norms we are going to consider.  The necessity to use the space $\G$ with the norm $|\cdot|$ is dictated by the needs of the heat kernel measure construction carried out in \cite{Gor00-1}, \cite{Gor00-2}, \cite{Gor02}, \cite{Gor03-2}.
\end{rem}

\begin{rem}\ll{rem:ad} Assumption \ref{E:norm} implies that the operator $\ad h$ is bounded on $\G$, namely, $\|\ad h\|\leqslant C|h|$ where the constant $C$ is as in Assumption \ref{E:norm}.
\end{rem}

\begin{thm}\ll{T:LocDiff}
If Assumption \ref{E:norm} is satisfied, then the exponential map is a diffeomorphism from a neighborhood of $0$ in $\G$ onto a neighborhood of $I$ in $\CM$.
\end{thm}

As before let $G_1 \subseteq G_2 \subseteq ... \subseteq G_n \subseteq ...\subseteq B(H)$ be a sequence of connected finite-dimensional Lie subgroups of $GL(H)$.

\begin{thm}\ll{T:INFTY=CM}

Suppose Assumption \ref{E:norm} holds. Then the group $\GG=\bigcup\limits_{n=1}^{\infty} G_n$ is dense in $\CM$ in the metric $d$, and the Cameron-Martin group $\CM$ is complete in the metric $d$. 

\end{thm}

In order to prove the surjectivity of the exponential map onto a neighborhood of the identity in the Cameron-Martin group $\CM$, it is necessary to prove that, for example, if $x$ and $y$ are small (in the norm $|\cdot|$) elements of $\G$, then $e^xe^y=e^z$ for some element $z$ in $\G$. This makes it unavoidable to use some version of the Baker-Campbell-Dynkin-Hausdorff (BCDH) formula (e.g. \cite{Dyn47}, \cite{Dyn49}).  

We begin by proving several preliminary results. Most of these lemmas were first proven in a somewhat different form in \cite{Gor00-2}. We give brief proofs here to make the exposition complete.  Lemma \ref{L:LOG} and Proposition \ref{L:LIP} are interesting in themselves. They show, with the help of the BCDH formula, that $\log$ has a derivative with values in $\G$. The standard integral formula for the logarithm $\log(I+x)=\int_0^1 x(I+sx)^{-1}ds$, although more useful in many contexts, does not easily yield values in $\G$. 

\begin{nota}\ll{log-def}
Let $g \in C^1_{CM}$ be such that  $\| g(t)-I\| <1$ for any $t \in [0,1]$.
Define the logarithm of $g$ by

\[
h(t)=\log g(t)=\sum\limits_{n=1}^{\infty} \frac{(-1)^{n-1}}{n}(g(t)-I)^n.
\] 
\end{nota}

\begin{prop}\ll{L:LIP}[The derivative of $log$ as a series]
Let $A(t)=g(t)^{-1}\Dot{g}(t)$, and for any $x \in \G$  define

\begin{equation} \ll{E:RIGHTSIDE OF ODE}
F(x,t)=A(t)+\frac 12[x, A(t)]-\sum\limits_{p=1}^{\infty}\frac{1}{2(p+2)p!}[...[[x,A(t)],
\overbrace{x],...,x}^{p}].
\end{equation}
Then the series converges in $\G$.
\end{prop}

\begin{proof}  Indeed, by Assumption \ref{E:norm}

\begin{multline*}
|F(x,t)|\leqslant 
|A(t)|+\frac{C|x| |A(t)|}{2}+
\sum\limits_{p=1}^{\infty}\frac{1}{2(p+2)p!}
\lnorm [...[[x,A(t)],\overbrace{x],...,x}^{p}]\rnorm
\leqslant
\\
|A(t)|+\frac{C|x| |A(t)|}{2}+
|A(t)|\sum\limits_{p=1}^{\infty}\frac{C^{p+1} |x|^{p+1}}{2(p+2)p!}<\infty, 
\end{multline*}
where $C$ is the same constant as in Assumption \ref{E:norm}. In particular, this means that $F(x, t) \in \G$ for any $x \in \G$ and any $0\leqslant t \leqslant 1$.
\end{proof}

\begin{lem}\ll{L:LOG}[The derivative of $log$]
Let $g \in C^1_{CM}$, $h=\log g$, then

\begin{equation}
\ll{E:ODE}
\Dot{h}(t)=F(h,t),  
\end{equation}
where $F(x,t)$ is defined by Equation \eqref{E:RIGHTSIDE OF ODE}.

\end{lem}

\begin{proof}
Indeed, $h(t+s)=\log g(t+s)=\log(g(t)g(t)^{-1}g(t+s))$. Let 
\[
f(t,s)=\log(g(t)^{-1}g(t+s)).
\] 
Then $h(t+s)=\log(e^{h(t)}e^{f(t,s)})=BCDH(h(t), f(t,s))$, where $BCDH(x,y)$ is given by the Baker-Campbell-Dynkin-Hausdorff formula for $x,y \in \g$

\begin{multline}
BCDH(x,y)=\log(\exp x \exp y)=
\\
\sum\limits_m \sum\limits_{p_i, q_i} 
\frac{(-1)^{m-1}}{m(\sum\limits_i p_i+q_i)}\frac{[...
\overbrace{[x,x],...,x}^{p_1}
],
\overbrace{y],...,y}^{q_1}
],
\overbrace{...,y],...y}^{q_m}
]}{p_1!q_1!...p_m!q_m!}.
\end{multline}
Note that

\begin{multline*}
\frac {df(t,s)}{ds}=
\lim_{\varepsilon \to 0}\frac{f(t,s+\varepsilon)-f(t,s)}{\varepsilon}=
\\
\sum\limits_{n=1}^{\infty} \frac{(-1)^{n-1}}{n} 
\sum\limits_{k=1}^{n}(g(t)^{-1}g(t+s)-I)^{k-1}g(t)^{-1}\Dot{g}(t+s)\left(g(t)^{-1}g(t+s)-I\right)^{n-k}, 
\end{multline*}
and therefore
 
\[
\frac {df(t,s)}{ds}\left|_{s=0}\right.=
g(t)^{-1}\Dot{g}=A(t).
\]

Note that $f(t,0)=0$, and therefore
\begin{multline*}
\Dot{h}(t)=\lim_{s \to 0}\frac{h(t+s)-h(t)}{s}=
\lim_{s \to 0}\frac{BCDH(h(t),f(t,s))-h(t)}{s}=
\\
A(t)+\frac 12[h(t), A(t)]+\frac{d}{ds}\left|_{s=0}\right. \sum\limits_{p} \frac{(-1)^{1}}{2(p+2)}
\frac{[...[[
h(t)
f(s,t)
],
\overbrace{h(t)],...,h(t)}^{p}
]}{1!1!p!0!}=
\\
A(t)+\frac 12[h(t), A(t)]-\sum\limits_{p=1}^{\infty}\frac{1}{2(p+2)p!}[...[[h(t),A(t)],
\overbrace{h(t)],...,h(t)}^{p}
]
\end{multline*}

\end{proof}

Denote by $d_n$ the distance metric on $G_n$ corresponding to the norm $|\cdot |$ restricted to the Lie algebra $\g _n$, and define the metric $d_{\infty}=\inf_n d_n$. 
As before let $\GG=\bigcup\limits_{n=1}^{\infty} G_n$. Then $\GG$ is a group contained in the Cameron-Martin group $\CM$. Moreover, for any $x,y \in \GG$ we have $d(x,y)\leqslant d_{\infty}(x,y)$.

\begin{lem}\ll{L:METRICS1}
Let $0<L<\ln 2/2C$, where $C$ is the same constant as in Assumption \ref{E:norm}. Then there is a positive constant $M$ such that

\begin{align*}
\left|d(I, e^x)-|x|\right| \leqslant  M |x|^2, \\
\left|d_{\infty}(I, e^x)-|x|\right| \leqslant  M |x|^2.
\end{align*}
for any $x \in \g$ provided $|x|<L$.

\end{lem}

\begin{proof}
First of all, the proof is the same for the metrics $d$ and $d_n$ with the same constants. We will show how to prove the estimate for the metric $d$. 
Joining $I$ to $e^x$ by the path $s \mapsto e^{sx}$, $0\leqslant s \leqslant 1$, we see that

\[
d(I,e^x)-|x|\leqslant \int_{0}^{1}|e^{-sx}\Dot{e^{sx}}|ds-|x|=0.
\]
Now let $g(s)$ and $h(s)=\log g(s)$ be as in Notation \ref{log-def}. As before, let $A(s)=g^{-1}(s)\Dot{g}(s)$, and therefore 
\[
A(s)=g^{-1}(s)\dot g(s)=\sum\limits_{k=0}^{\infty}\frac{(-\ad  h)^k(\dot h)}{(k+1)!}=
\dot h-[h,\dot h]+\frac{[h,[h,\dot h]]}{2!}-\frac{[h,[h,[h,\dot h]]]}{3!}+...
\]
Then by Remark \ref{rem:ad} for any smooth path $g(s)$ from $I$ to $e^x$

\begin{multline*}
0\leqslant |x|-d(I,e^x)=\left|\int_{0}^{1}\dot h(s)ds\right|-\int_{0}^{1} \left|g^{-1}(s)\dot g(s)\right|ds \leqslant \\
\left|\int_{0}^{1}\dot h(s)ds\right|-\left|\int_{0}^{1} g^{-1}(s)\dot g(s)ds\right|
\leqslant \left|\int_{0}^{1}\dot h(s)-g^{-1}(s)\dot g(s)ds\right|=
\\
\left|\int_{0}^{1}\sum\limits_{k=1}^{\infty}\frac{(-\ad  h)^k(\dot h)}{(k+1)!}ds\right|\leqslant \sum\limits_{k=1}^{\infty}\frac{C^k}{(k+1)!}\int_{0}^{1}|h(s)|^k|\dot h(s)|ds
,
\end{multline*}
where the constant $C$ is the same as in Assumption \ref{E:norm}. In particular, if $h(s)=sx, \ 0\leqslant s \leqslant 1$, we have

\[
\lnorm d(I,e^x)-|x|\rnorm \leqslant \sum\limits_{k=1}^{\infty}\frac{C^k |x|^{k+1}}{(k+1)!}\int_{0}^{1}s^k ds \leqslant |x|^2 \sum\limits_{k=1}^{\infty}\frac{C^k L^{k-1}}{(k+1)!(k+1)}.
\]
Thus we can take

\[
M=\sum\limits_{k=1}^{\infty}\frac{C^k L^{k-1}}{(k+1)!(k+1)}.
\]
Note that the above proof would go through for $d_n$, where $n$ is such that $x \in \g_n$. Thus the statement of Lemma \ref{L:METRICS} holds for $d$ replaced by $d_{\infty}$ with the same constant $M$.
\end{proof}

\begin{cor}\ll{L:METRICS}
By choosing $h(s)=sy+(1-s)x$ for $0 \leqslant s \leqslant 1$ in the above proof we have that

\begin{align*}
\left|d(e^x, e^y)-|x-y|\right| \leqslant  M |x-y|^2, \\
\left|d_{\infty}(e^x, e^y)-|x-y|\right| \leqslant  M |x-y|^2.
\end{align*}
for any $x, y \in \g$ provided $|x|<L$ and $|y|<L$.

\end{cor}

\begin{lem}\ll{L:LOCAL}
Let $g \in C^1_{CM}$. Suppose that 

\[
g(0)=I, \ \| g(t)-I\|<1 \mbox{ \ for \ any \ } t \in [0,1]  \mbox{ \ and \ }  |\log g(1)|<L,
\]
where $L$ is the same as in Corollary \ref{L:METRICS}. Then there is a sequence $g_n \in G_n$ such that $\lim\limits_{n\to \infty}d(g_n, g(1))=0$.
\end{lem}

\begin{proof}
As before let $h(s)=\log g(s)$. 
Note that by Proposition \ref{L:LIP} $h(t) \in \G$. Therefore there are $h_n \in \g_n$ such that 
$|h(1)-h_n|\xrightarrow[n \to \infty]{} 0$. Let $g_n=e^{h_n}$. 
By Corollary \ref{L:METRICS} we have
\[
d(e^{h(1)},g_n)\leqslant C_2 |h(1)-h_n|\xrightarrow[n \to \infty]{} 0.
\]
By a direct calculation we have $e^{h(t)}=g(t)$. 
\end{proof}

\begin{prop}\ll{R:METRICS} 
Suppose Assumption \ref{E:norm} holds, then 

\[
d(x,y)=d_{\infty}(x,y)
\]
for any $x, y \in \GG$.
\end{prop}

\begin{proof} Both metrics $d$ and $d_{\infty}$ are left-invariant on $\GG$, therefore we can assume that $x=I$. As we pointed out earlier, for any $y \in \GG$ we have $d(I,y)\leqslant d_{\infty}(I,y)$. Therefore we only need to prove that $d_{\infty}(I,y)\leqslant d(I,y)$ for any $y \in \GG$. 

Denote $d(I,y)=D$. For any $\e'>0$ there is a path $g(s) \in C^1_{CM}$ such that $g(0)=I, g(1)=y$ and

\[
D \leqslant \int_0^1 |g^{-1}(s)\dot g(s)|ds\leqslant D+\e'.
\]
Denote $\e=\int_0^1 |g^{-1}(s)\dot g(s)|ds-D$. Let $m>1, m \in \mathbb N$. There exist $ t_k,  0=t_0 < t_1 <  ...t_{m-1} < t_m=1$, $k=0,1,...,m$ such that 

\[
\int_{t_k}^{t_{k+1}} |g^{-1}(s)\dot g(s)|ds=\frac{D+\e}{m} \mbox{ \ and \ }|\log y_{k+1}^{-1}y_k|\leqslant L,
\]
where $y_k=g(t_k)$ and $L<\min\{\ln 2/2C, 1/M\}$. Note that in general $y_k$ is in the Cameron-Martin group $\CM$, not in the group $\GG$. By Lemma \ref{L:LOCAL} there exist $x_k \in \GG$ such that $d(y_k,x_k)<\e/m$. By applying Corollary \ref{L:METRICS} twice we have that

\begin{multline*}
d_{\infty}(x_k,x_{k+1})\leqslant  |\log x_{k+1}^{-1}x_k|+M|\log x_{k+1}^{-1}x_k|^2 \leqslant 
\\d(x_k,x_{k+1})+2M|\log x_{k+1}^{-1}x_k|^2
\leqslant
d(x_k,x_{k+1})+2M|\log y_{k+1}^{-1}y_k|^2
.
\end{multline*}
Corollary \ref{L:METRICS} also implies that

\[
|\log y_{k+1}^{-1}y_k|-M|\log y_{k+1}^{-1}y_k|^2\leqslant d(y_k, y_{k+1}) \leqslant \frac{D+\e}{m},
\]
and since $L<1/M$

\[
|\log y_{k+1}^{-1}y_k| \leqslant \frac{D+\e}{m(1-ML)}.
\]
Thus

\begin{multline*}
d_{\infty}(x_k,x_{k+1})\leqslant  
d(x_k,x_{k+1})+2M\left(\frac{D+\e}{m}\right)^2
\leqslant \\
d(x_k,y_k)+d(y_k,y_{k+1})+d(y_{k+1},x_{k+1})+2M\left(\frac{D+\e}{m}\right)^2
\leqslant 
\\
d(y_k,y_{k+1})+\frac{2\e}{m}+2M\left(\frac{D+\e}{m}\right)^2
\leqslant
\frac{D+\e}{m}+\frac{2\e}{m}+2M\left(\frac{D+\e}{m}\right)^2.
\end{multline*}
Finally,

\[
d_{\infty}(I,y)\leqslant m\left(\frac{D+\e}{m}+\frac{2\e}{m}+2M\left(\frac{D+\e}{m}\right)^2\right)=
D+3\e+\frac{2M(D+\e)^2}{m}.
\]
Recall that $D=d(I,y)$, and $\e$ and $m$ are arbitrary.
\end{proof}

Theorem \ref{T:LocDiff} now is a direct consequence of Lemma \ref{L:LOG} and Lemma \ref{L:METRICS} which imply that the exponential and logarithmic functions are well-defined  and differentiable in neighborhoods of the identity and zero respectively.

\begin{proof}[\textbf{Proof of Theorem \ref{T:INFTY=CM}}]

First of all, $\GG \subseteq \CM$ since $G_n \subseteq \CM$ for all $n$. Let $g \in \GG$, $k \in \CM$, and suppose we have a path $g(t):[0,1] \to \CM, g(t) \in C^1_{CM}, g(0)=g, g(1)=k$. Without a loss of generality we can assume that $\|g(t)-g\|<1$, $d(\log k,\log g)<L$, where $L$ is the same as in Lemma \ref{L:METRICS}. Otherwise the path $\log g(t)$ can be divided into a finite number of subpaths satisfying the condition.  Lemma \ref{L:LOCAL} implies that for any $\e >0$ there is $m \in \GG$ such that $d(m, g^{-1}k)<\e$. Then

\[
d(m, g^{-1}k)=d(gm, k) <\e.
\]
Note that $gm \in \GG$, and therefore we have shown that elements of $\CM$ can be approximated by elements of $\GG$. 
\end{proof}

\section{Is $\G$ a Lie algebra?}\ll{S:Lie}
 
As in the previous section, let $H$ be a separable Hilbert space, $B(H)$ be the space of bounded operators on $H$, and $I$ be the identity operator. Here we restrict ourselves to the case of the Hilbert-Schmidt groups. By $HS$ we denote the space of Hilbert-Schmidt operators on $H$. The space $HS$ is equipped with the Hilbert-Schmidt inner product $(\cdot, \cdot )_{HS}$. Let $GL(H)$ be the group of invertible elements of $B(H)$. Then a Hilbert-Schmidt group is a closed subgroup of $GL(H)$ such that $A-I \in HS$ for any $A\in G$. Note that the set $\{A\in GL(H): A-I \in HS\}$ is a group. 

In Section \ref{S:CM} we considered a sequence of finite-dimensional groups $G_n$. Suppose now that in addition $G_n \subset I+HS$. Then their Lie algebras satisfy $\g_n \subset HS$. As before we assume that the infinite-dimensional Lie algebra $\g=\bigcup\limits_{n=1}^{\infty} \g_n$ is equipped with a Hermitian inner product $(\cdot,\cdot)$, and the corresponding norm is denoted by $|\cdot |$. 
Throughout this section we assume the following modified version of Assumption \ref{A:compl}.
\begin{as}[Hilbert-Schmidt Completeness Assumption]\ll{A:complHS}
There is a subspace $\G$ of $HS$ such that the Lie algebra $\g$ is contained in $\G$ and the given inner product $(\cdot,\cdot)$ on $\g$ extends to $\G$, which is complete with respect to this inner product. As before we will abuse notation by using $(\cdot,\cdot)$ to denote the extended inner product on $\G$ and by $|\cdot |$ the corresponding norm. We assume that $\g$ is dense in $\G$ in the norm $|\cdot |$. 
\end{as}

\begin{rem}\ll{R:Q} In our earlier papers we assumed that $(\cdot, \cdot)_{\G}$ is given by $(x, y)_{\G}=(x, Q^{-1}y)_{HS}$, where $Q$ is a one-to-one nonnegative trace class operator  on $HS$ for which each $\g_n$ is an invariant subspace. The assumption that $Q$ is trace-class assures that the heat kernel measure constructed in our previous work (\cite{Gor00-1}, \cite{Gor00-2}, \cite{Gor02}, \cite{Gor03-2}) actually lives in $HS+I$. In the present paper we do not assume that $Q$ is trace-class unless it is stated explicitly. Moreover, we do not use the operator $Q$, but rather describe the assumptions on $Q$ in terms of an orthonormal basis of $\G$. 
\end{rem}
In the next statement we use the fact that we can view an element of $HS$ as an infinite matrix $A=\{a_{ij}\}_{i,j=1}^{\infty}$ such that the sum $\sum_{i,j}|a_{ij}|^2$ is finite. 
Then $e_{ij}$, the matrices with $1$ at the $ij$th place and $0$ at all other places, form an orthonormal basis of $HS$ with the inner product $( \cdot, \cdot )_{HS}$. Let us describe an example of the setting introduced above. Namely, let $\G$ be the vector space generated by the orthonormal basis $\xi_{ij}=\l_{ij}e_{ij}, \ (i,j)\in A \subseteq \mathbb N \times \mathbb N$ for some $\l_{ij}>0$. Then the inner product on $\G$ is determined by
$(\xi_{ij}, \xi_{km})_{\G}=(e_{ij}, e_{km})_{HS}$. It turns out that $\G$ might not be a Lie algebra.

\begin{prop} 
There exists a sequence of positive numbers $\l_{ij}$ such that $\G$ is not a Lie algebra.
\end{prop}

\begin{proof}

Let $\l_{ij}=\l_{ji}$ for any $i, j$, and $x, y \in \G$ be such that

\[
x=\sum_{l=0}^{\infty} x_l (\xi_{3l+1,3l+2}-\xi_{3l+2,3l+1}), \
y=\sum_{l=0}^{\infty} y_l (\xi_{3l+1,3l+3}-\xi_{3l+3,3l+1}).  
\]
Then

\begin{align*}
[x,y]&=\sum_{l,n} x_l y_n [\xi_{3l+1,3l+2}-\xi_{3l+2,3l+1}, \xi_{3n+1,3n+3}-\xi_{3n+3,3n+1}]
\\
&=-\sum_{l} x_l y_l \z{\l_{3l+1,3l+2}\l_{3l+1,3l+3}}{\l_{3l+2,3l+3}}(\xi_{3l+2,3l+3}-\xi_{3l+3,3l+2}).
\end{align*}

Thus 

\[
|[x,y]|^2=\sum_{l} x_l^2 y_l^2
\z{\l_{3l+1,3l+2}^2\l_{3l+1,3l+3}^2}{\l_{3l+2,3l+3}^2}.
\]

Let $x_l=y_l=\l_{3l+1,3l+2}=\l_{3l+1,3l+3}=a_l$ and $\l_{3l+2,3l+3}=b_l$, 
where $a_l$ and $b_l$ are $\ell^2$-sequences. Then 
$|[x,y]|^2=\sum_{l} a_l^8 b_l^{-2}=\infty$ if, for example, $a_l=1/l$ and $b_l=1/l^4$. 

\end{proof}

The next result shows that there exist inner products such that Continuity Assumption \ref{E:norm} on the Lie bracket is satisfied.
 
\begin{thm}\ll{T:split} Suppose $\l_{i,j}=\l_{i}\l_{j}$, $i, j \in \mathbb{N}$, then for any $x,y \in \G$

\[
|[x,y]|\leqslant 2\sup_i{\l_i^2}\  |x||y|.
\]
\end{thm}

\begin{proof}
Let $x=\sum_{i,j=1}^{\infty} x_{i,j}\xi_{ij}$ and $y=\sum_{k,m=1}^{\infty} y_{k,m}\xi_{km}$. Then

\[
xy=\sum_{i,j,m} x_{i,j}y_{j,m} \z{\l_{i,j}\l_{j,m}}{\l_{i,m}}\xi_{im}
\]
and therefore

\begin{align*}
|xy|^2&=\sum_{i,m} \left(\sum_{j}x_{i,j}y_{j,m} \z{\l_{i,j}\l_{j,m}}{\l_{i,m}}\right)^2
\\
&=\sum_{i,m} \left(\sum_{j}x_{i,j}y_{j,m} \l_{j}^2\right)^2 \leqslant  \sum_{i,j}x_{i,j}^2\l_{j}^2 \sum_{k,m}y_{k,m}^2 \l_{k}^2.
\end{align*}

\end{proof}

\begin{cor}
If $\l_{i,j}=\l_{i}\l_{j}$, $i, j \in \mathbb{N}$, and $\{\l_{i}\}_{i=1}^{\infty}   \in \ell^p$ for any $p>0$, then Assumption \ref{E:norm} on the Lie bracket is satisfied. In particular, if $p=2$, then the operator $Q$ mentioned in Remark \ref{R:Q} is trace-class. 
\end{cor}

\section{Riemannian geometry of the Hilbert-Schmidt groups: definitions and preliminaries}\ll{S:geometry}

The goal of the next two sections is to see if there exists a natural Lie algebra for $\CM$ and an inner product on it such that the Ricci curvature is bounded from below. The first obstacle in answering such a question is the absence of geometric definitions. We chose to follow the work of J. Milnor for finite-dimensional Lie groups in \cite{Milnor76}. There he described the Riemannian geometry of a Lie group with a Riemannian metric invariant under left translation. One of his aims was to see how the choice of an orthonormal basis of the (finite-dimensional) Lie algebra determines the curvature properties of the corresponding Lie group. This is exactly the question we study, but in infinite dimensions: how the choice of the inner product on $\g$ changes the Riemannian geometry of the group $\CM$. We consider general norms on $\g$ which are diagonal in a certain sense. This allows us to compute the Ricci curvature in two important cases: the first case is when the norm on $\g$ is the Hilbert-Schmidt norm, and the second one is when the norm on $\g$ is determined by a nonnegative trace class operator on $HS$. The latter assumption assures that the corresponding heat kernel measure constructed in our previous work (\cite{Gor00-1}, \cite{Gor00-2}, \cite{Gor02}, \cite{Gor03-2}) actually lives in $HS+I$. We use finite-dimensional approximations to $\g$ to define the sectional and Ricci curvatures. Our results show that for the general, orthogonal and upper triangular Hilbert-Schmidt algebras the Ricci curvature generally is not bounded from below. Moreover, for the upper triangular Hilbert-Schmidt algebra the Ricci curvature is identically minus infinity. 

Let $\g$ be a infinite-dimensional Lie algebra equipped with an inner product $(\cdot ,\cdot )$. We assume that $\g$ is complete. By Theorem \ref{T:split} this is the case for all examples we consider later in this Section.

\begin{df}
The \textbf{Levi-Civita connection} $\n_x$ is defined by
\[
(\n_x y,z)=\z 12(([x, y],z)-([y, z],x)+([z, x],y))
\]
for any $x,y,z \in \g$.
\end{df}

\begin{df}
\begin{enumerate}
\item The \textbf{Riemannian curvature tensor} $R$ is defined by
\[
R_{xy}=\n_{[x,y]}-\n_x\n_y+\n_y\n_x, \ x,y \in \g.
\]

\item For any orthogonal $x,y$ in $\g$
\[
K(x,y)=(R_{xy}(x),y)
\]
is called the \textbf{sectional curvature}.

\item Let $\{\xi_i\}_{i=1}^{\infty}$ be an orthonormal basis of $\g$, $N$ be finite, then

\[
R^N(x)=\sum_{i=1}^{N}K(x,\xi_i)=\sum_{i=1}^{N}(R_{x\xi_i}(x),\xi_i)
\]
is the \textbf{truncated Ricci curvature}.

\item

Let $N$ be finite, then 

\[
\hat R^N(x)=\sum_{i=1}^{N}R_{\xi_i x}(\xi_i)
\]
is the \textbf{truncated self-adjoint Ricci curvature or transformation}.

\end{enumerate}
\end{df}
The self-adjoint Ricci transformation is a convenient computational tool. First of all, 

\[
R^N(x)=(\hat R^N(x),x).
\]
Then if $\{\xi_i\}_{i=1}^{\dim \g}$ is an orthonormal basis which diagonalizes $\hat R$, that is,  $\hat R^N(\xi_i)=a_i \xi_i$, then 

\[
R^N(x)=\sum_{i=1}^{N} a_i x_i^2, \ x=\sum_{i=1}^{dim \g} x_i \xi_i.
\]
The numbers $a_i$ are called the principal Ricci curvatures.

As in Section \ref{S:Lie}, let $\{e_{ij}\}_{i,j=1}^{\infty}$ be the standard basis of the space of Hilbert-Schmidt operators $HS$. The Lie bracket for these basis elements can be written as

\[
[e_{ij},e_{km}]=\y_{j,k}e_{im}-\y_{i,m}e_{kj}
\]
where $\y_{ln}$ is Kronecker's symbol. As before, we study a subspace $\G$ of $HS$ generated by an orthonormal basis $\xi_{ij}=\l_i\l_je_{ij}$ for some $\l_i>0$, \ $(i,j) \in A \subseteq \mathbb N \times \mathbb N$.

\section{Ricci curvature}\ll{S:Ricci}

The main results of this section are Theorems \ref{T:GL}, \ref{T:O} and \ref{T:T}. For the skew-symmetric and triangular infinite matrices, the orthonormal basis we use actually diagonalizes the truncated self-adjoint Ricci curvature. Then the results of this section show that if we define the Ricci curvature as the limit of the truncated Ricci curvature as the dimension goes to $\infty$, it is not bounded from below. Moreover, for the upper triangular matrices the Ricci curvature is identically negative infinity.

\subsection{General Hilbert-Schmidt algebra} 

Let $\{\l_i\}_{i=1}^{\infty}$ be a bounded sequence of strictly positive numbers. 
In this section we consider the infinite dimensional Lie algebra $\G$ generated by the orthonormal basis $\xi_{ij}=\l_i\l_je_{ij}$. Then $\G$ is a Lie subalgebra of $\mathfrak{g}\mathfrak{l}_{HS}$. Recall that if the sequence $\{\l_i\}_{i=1}^{\infty}$ is bounded, then by Theorem \ref{T:split}, Continuity Assumption \ref{E:norm} is satisfied for the corresponding  norm on $\G$.

\begin{thm}\ll{T:GL}
\begin{enumerate}

\item
Let $N>\max{i,j}$, then the truncated Ricci curvature is

\[
R^N_{ij}
=
R^N(\xi_{ij})
=
\frac 14
(
6\y_{i,j}\l_{i}^4        
-4\y_{i,j}\l_{i}^4 N     
-2\l_{i}^4 N    
-2\l_{j}^4  N   
+2\sum_{m=1}^{N} \l_{m}^4    
)
.
\]

\item

Suppose that $\{\l_i\}_{i=1}^{\infty} \in \ell^2$. Then 
\[
\lim_{N \to \infty} R^N(\xi_{ij})
=
\frac 12
\lim_{N \to \infty} (    
-\l_{i}^4 N    
-\l_{j}^4  N   
+\sum_{m=1}^{N} \l_{m}^4    
)=-\infty
,
\]
if $i\not=j$.

\[
\lim_{N \to \infty} R^N(\xi_{ii})
=
\frac 12
\lim_{N \to \infty}(
3\l_{i}^4        
-4\l_{i}^4 N       
+\sum_{m=1}^{N} \l_{m}^4    
)=-\infty
.
\]

\item
For the Hilbert-Schmidt inner product, the truncated Ricci curvature is

\[
R^N_{ij}
=
\frac 12
(
3\y_{i,j}        
-2\y_{i,j} N     
-N        
)
.
\]
\end{enumerate}

\end{thm}

This theorem is a direct consequence of the following results. First of all, note that the Lie bracket can be written

\[
[\xi_{ij},\xi_{km}]=\y_{j,k}\l_{j}^2\xi_{im}-
\y_{i,m}\l_{i}^2\xi_{kj}, 
\]
where $\y_{l,n}$ is Kronecker's symbol. Denote $\n_{ab}=\n_{\xi_{ab}}$, then

\begin{lem}

\[
\n_{ab}\xi_{cd}=
\frac 12(
\y_{b,c}\l_{b}^2\xi_{ad}
-\y_{a,d}\l_{a}^2\xi_{cb}
-\y_{a,c}\l_{d}^2\xi_{db}
+\y_{b,d}\l_{c}^2\xi_{ac}
+\y_{b,d}\l_{a}^2\xi_{ca}
-\y_{a,c}\l_{b}^2\xi_{bd}
)
\]

\end{lem}

\begin{proof}
\begin{multline*}
(\n_{ab}\xi_{cd},\xi_{ef})=\frac 12\left(([\xi_{ab},\xi_{cd}],\xi_{ef})-([\xi_{cd},\xi_{ef}],\xi_{ab})+([\xi_{ef},\xi_{ab}],\xi_{cd})\right)=\\
\frac 12(
\y_{b,c}\l_{b}^2(\xi_{ad},\xi_{ef})
-\y_{a,d}\l_{a}^2(\xi_{cb},\xi_{ef})
-\y_{d,e}\l_{d}^2(\xi_{cf},\xi_{ab})
\\
+\y_{c,f}\l_{c}^2(\xi_{ed},\xi_{ab})
+\y_{f,a}\l_{f}^2(\xi_{eb},\xi_{cd})
-\y_{e,b}\l_{e}^2(\xi_{af},\xi_{cd})
)
=\\
\frac 12(
\y_{b,c}\y_{e,a}\y_{f,d}\l_{b}^2
-\y_{a,d}\y_{e,c}\y_{f,b}\l_{a}^2
-\y_{a,c}\y_{e,d}\y_{f,b}\l_{d}^2
\\
+\y_{b,d}\y_{e,a}\y_{f,c}\l_{c}^2
+\y_{b,d}\y_{e,c}\y_{f,a}\l_{a}^2
-\y_{a,c}\y_{f,d}\y_{e,b}\l_{b}^2
).
\end{multline*}

\end{proof}

\begin{lem}
\begin{multline*}
R_{ij,km}\xi_{ij}
=\\
\frac 14
(
3\y_{j,k}\y_{j,m}\l_{i}^2\l_{j}^2\xi_{ii}
+2\y_{i,k}\y_{j,m}(\l_{i}^4+\l_{j}^4)\xi_{ij}
+2\y_{i,m}\y_{j,k}\l_{i}^2\l_{j}^2\xi_{ij}
\\
+2\y_{i,k}\y_{j,m}\l_{i}^2\l_{j}^2\xi_{ji}
+3\y_{i,k}\y_{i,m}\l_{i}^2\l_{j}^2\xi_{jj}
\\
-\y_{i,m}\l_{i}^2\l_{k}^2\xi_{ik}
-2\y_{i,j}\y_{i,m}\l_{i}^2\l_{k}^2\xi_{ik}
+\y_{i,m}\l_{i}^4\xi_{ki}
-2\y_{i,j}\y_{i,m}\l_{i}^4\xi_{ki}
\\
-4\y_{i,k}\l_{i}^4\xi_{im}
+\y_{i,k}\l_{j}^4\xi_{im}
+\y_{i,k}\l_{m}^4\xi_{im}
-2\y_{i,j}\y_{i,k}\l_{i}^4\xi_{im}
-\y_{i,k}\l_{i}^2\l_{m}^2\xi_{mi}
-2\y_{i,j}\y_{i,k}\l_{i}^2\l_{m}^2\xi_{mi}
\\
-\y_{j,m}\l_{j}^2\l_{k}^2\xi_{jk}
+\y_{j,m}\l_{i}^4\xi_{kj}
+\y_{j,m}\l_{k}^4\xi_{kj}
-4\y_{j,m}\l_{j}^4\xi_{kj}
+\y_{j,k}\l_{j}^4\xi_{jm}
-\y_{j,k}\l_{j}^2\l_{m}^2\xi_{mj}
)
\end{multline*}

\end{lem}

\begin{proof}
The Riemannian curvature tensor applied to $\xi_{ij}$ is

\begin{multline*}
R_{\xi_{ij} \xi_{km}}\xi_{ij}=R_{ij,km}\xi_{ij}
=
\n_{[\xi_{ij},\xi_{km}]}\xi_{ij}-\n_{ij}\n_{km}\xi_{ij}+
\n_{km}\n_{ij}\xi_{ij}
=\\
\y_{j,k}\l_{j}^2\n_{im}\xi_{ij}
-\y_{i,m}\l_{i}^2\n_{kj}\xi_{ij}
-\n_{ij}\n_{km}\xi_{ij}+
\n_{km}\n_{ij}\xi_{ij}
=\\
\frac 12\y_{j,k}\l_{j}^2
(
\y_{i,m}\l_{i}^2\xi_{ij}
-\y_{i,j}\l_{i}^2\xi_{im}
-\l_{j}^2\xi_{jm}
-\l_{m}^2\xi_{mj}
+2\y_{m,j}\l_{i}^2\xi_{ii}
)\\
-\frac 12\y_{i,m}\l_{i}^2
(
\l_{i}^2\xi_{ki}
+\l_{k}^2\xi_{ik}
+\y_{j,i}\l_{j}^2\xi_{kj}
-\y_{k,j}\l_{k}^2\xi_{ij}
-2\y_{i,k}\l_{j}^2\xi_{jj}
)\\
-\frac 12\n_{ij}
(
\y_{i,m}\l_{m}^2\xi_{kj}
-\y_{k,j}\l_{k}^2\xi_{im}
-\y_{i,k}\l_{j}^2\xi_{jm}
+\y_{m,j}\l_{i}^2\xi_{ki}
+\y_{m,j}\l_{k}^2\xi_{ik}
-\y_{i,k}\l_{m}^2\xi_{mj}
)
\\
+\frac 12\n_{km}
(
\l_{i}^2\xi_{ii}
+\l_{i}^2\xi_{ii}
+\y_{i,j}\l_{j}^2\xi_{ij}
-\y_{i,j}\l_{i}^2\xi_{ij}
-2\l_{j}^2\xi_{jj}
)
\end{multline*}

\end{proof}

\begin{lem}
The sectional curvature is

\begin{multline*}
K(\xi_{ij},\xi_{km})=(R_{ij,km}\xi_{ij},\xi_{km})
=\\
\frac 14
(
6\y_{i,k}\y_{i,m}\y_{j,k}\y_{j,m}\l_{i}^4       
+2\y_{i,k}\y_{j,m}\l_{i}^4     
+2\y_{i,k}\y_{j,m}\l_{j}^4     
-2\y_{i,k}\y_{i,m}\y_{k,m}\l_{i}^4      
\\  
+\y_{i,m}\l_{i}^4     
-2\y_{i,j}\y_{i,m}\l_{i}^4     
-4\y_{i,k}\l_{i}^4     
+\y_{i,k}\l_{j}^4     
+\y_{i,k}\l_{m}^4    
-2\y_{i,j}\y_{i,k}\l_{i}^4     
-2\y_{j,k}\y_{j,m}\y_{m,k}\l_{j}^4
\\
+\y_{j,m}\l_{i}^4    
+\y_{j,m}\l_{k}^4    
-4\y_{j,m}\l_{j}^4     
+\y_{j,k}\l_{j}^4     
).
\end{multline*}
\end{lem}

\begin{proof}
\begin{multline*}
(R_{ij,km}\xi_{ij},\xi_{km})
=\\
\frac 14
(
3\y_{i,k}\y_{i,m}\y_{j,k}\y_{j,m}\l_{i}^2\l_{j}^2     
+2\y_{i,k}\y_{j,m}\l_{i}^4     
+2\y_{i,k}\y_{j,m}\l_{j}^4     
\\
+2\y_{i,k}\y_{i,m}\y_{j,k}\y_{j,m}\l_{i}^2\l_{j}^2     
+2\y_{i,k}\y_{i,m}\y_{j,k}\y_{j,m}\l_{i}^2\l_{j}^2     
+3\y_{i,k}\y_{i,m}\y_{j,k}\y_{j,m}\l_{i}^2\l_{j}^2     
\\
-\y_{i,k}\y_{i,m}\y_{k,m}\l_{i}^2\l_{k}^2     
-2\y_{i,j}\y_{i,k}\y_{i,m}\y_{k,m}\l_{i}^2\l_{k}^2     
+\y_{i,m}\l_{i}^4     
-2\y_{i,j}\y_{i,m}\l_{i}^4     
\\
-4\y_{i,k}\l_{i}^4     
+\y_{i,k}\l_{j}^4     
+\y_{i,k}\l_{m}^4    
-2\y_{i,j}\y_{i,k}\l_{i}^4     
-\y_{i,k}\y_{i,m}\y_{m,k}\l_{i}^2\l_{m}^2     
-2\y_{i,j}\y_{i,k}\y_{i,m}\y_{m,k}\l_{i}^2\l_{m}^2    
\\
-\y_{j,k}\y_{j,m}\y_{k,m}\l_{j}^2\l_{k}^2    
+\y_{j,m}\l_{i}^4    
+\y_{j,m}\l_{k}^4    
-4\y_{j,m}\l_{j}^4     
+\y_{j,k}\l_{j}^4    
-\y_{j,k}\y_{j,m}\y_{m,k}\l_{j}^2\l_{m}^2     
)
=\\
\frac 14
(
6\y_{i,k}\y_{i,m}\y_{j,k}\y_{j,m}\l_{i}^4       
+2\y_{i,k}\y_{j,m}\l_{i}^4     
+2\y_{i,k}\y_{j,m}\l_{j}^4     
-2\y_{i,k}\y_{i,m}\y_{k,m}\l_{i}^4      
\\  
+\y_{i,m}\l_{i}^4     
-2\y_{i,j}\y_{i,m}\l_{i}^4     
-4\y_{i,k}\l_{i}^4     
+\y_{i,k}\l_{j}^4     
+\y_{i,k}\l_{m}^4    
-2\y_{i,j}\y_{i,k}\l_{i}^4     
-2\y_{j,k}\y_{j,m}\y_{m,k}\l_{j}^4
\\
+\y_{j,m}\l_{i}^4    
+\y_{j,m}\l_{k}^4    
-4\y_{j,m}\l_{j}^4     
+\y_{j,k}\l_{j}^4     
).
\end{multline*}

\end{proof}

\subsection{Orthogonal Hilbert-Schmidt algebra} Note that $b_{ij}=(e_{ij}-e_{ji})/\sqrt{2}, i<j$ is an orthonormal basis for the space of skew-symmetric Hilbert-Schmidt operators $\mathfrak{s}\mathfrak{o}_{HS}$. Let $\{\l_i\}_{i=1}^{\infty}$ be a bounded sequence of strictly positive numbers. In this section we consider the infinite dimensional Lie algebra $\G$ generated by the orthonormal basis $\xi_{ij}=\l_i\l_jb_{ij}$. Recall that if the sequence $\{\l_i\}_{i=1}^{\infty}$ is bounded, then by Theorem \ref{T:split}, Continuity Assumption \ref{E:norm} is satisfied for the corresponding  norm on $\G$.

In what follows the convention is that $\xi_{ij}=0$, if $i\geqslant j$. The Lie bracket for these basis elements can be written as

\[
[\xi_{ij},\xi_{km}]=\frac{1}{\sqrt{2}}(
\y_{j,k}\l_{j}^2\xi_{im}
+\y_{j,m}\l_{j}^2(\xi_{ki}-\xi_{ik})
+\y_{i,k}\l_{i}^2(\xi_{mj}-\xi_{jm})
-\y_{i,m}\l_{i}^2\xi_{kj}
)
\]
We begin with several computational lemmas.

\begin{lem}

\begin{multline*}
\n_{ab}\xi_{cd}=\\
\frac{1}{2\sqrt{2}}(
\l_{b}^2\y_{b,c}\xi_{a,d}
+ \l_{b}^2\y_{b,d}(\xi_{c,a}-\xi_{a,c})
+ \l_{a}^2\y_{a,c}(\xi_{d,b}-\xi_{b,d})  
-\l_{a}^2\y_{a,d}\xi_{c,b}
\\
-\l_{d}^2\y_{a,c}\xi_{d,b}
-\l_{d}^2(\y_{b,c}\xi_{a,d}-\y_{a,c}\xi_{b,d})  
-\l_{c}^2(\y_{b,d}\xi_{c,a}-\y_{a,d}\xi_{c,b})  
+\l_{c}^2\y_{b,d}\xi_{a,c}
\\
+\l_{a}^2\y_{b,d}\xi_{c,a}
+\l_{b}^2(\y_{a,c}\xi_{d,b}-\y_{a,d}\xi_{c,b}) 
+ \l_{a}^2(\y_{b,c}\xi_{a,d}-\y_{b,d}\xi_{a,c})  
-\l_{b}^2\y_{a,c}\xi_{b,d}).
\end{multline*}

\end{lem}

\begin{proof}
\begin{multline*}
(\n_{ab}\xi_{cd},\xi_{ef})
=\\
\frac {1}{2\sqrt{2}}(
((
\y_{b,c}\l_{b}^2\xi_{ad}
+\y_{b,d}\l_{b}^2(\xi_{ca}-\xi_{ac})
+\y_{a,c}\l_{a}^2(\xi_{db}-\xi_{bd})
-\y_{a,d}\l_{a}^2\xi_{cb}
)
,\xi_{ef})\\
-(
\y_{d,e}\l_{d}^2\xi_{cf}
+\y_{d,f}\l_{d}^2(\xi_{ec}-\xi_{ce})
+\y_{c,e}\l_{c}^2(\xi_{fd}-\xi_{df})
-\y_{c,f}\l_{c}^2\xi_{ed}
,\xi_{ab})\\
+(
\y_{f,a}\l_{f}^2\xi_{eb}
+\y_{f,b}\l_{f}^2(\xi_{ae}-\xi_{ea})
+\y_{e,a}\l_{e}^2(\xi_{bf}-\xi_{fb})
-\y_{e,b}\l_{e}^2\xi_{af}
,\xi_{cd})
)=\\
\frac{1}{2\sqrt{2}}(
\l_{b}^2\y_{a,e}\y_{b,c}\y_{d,f}
+ \l_{b}^2\y_{b,d}(\y_{c,e}\y_{a,f}-\y_{a,e}\y_{c,f})
+ \l_{a}^2\y_{a,c}(\y_{d,e}\y_{b,f}-\y_{d,f}\y_{b,e})  
\\
-\l_{a}^2\y_{a,d}\y_{b,f}\y_{c,e}
-\l_{d}^2\y_{a,c}\y_{b,f}\y_{d,e}
-\l_{d}^2\y_{d,f}(\y_{a,e}\y_{b,c}-\y_{a,c}\y_{b,e})  
-\l_{c}^2\y_{c,e}(\y_{a,f}\y_{b,d}-\y_{a,d}\y_{b,f})  
\\
+\l_{c}^2\y_{a,e}\y_{b,d}\y_{c,f}
+\l_{a}^2\y_{a,f}\y_{b,d}\y_{c,e}
+\l_{b}^2\y_{b,f} (\y_{a,c}\y_{e,d}-\y_{e,c}\y_{a,d}) 
+ \l_{a}^2\y_{a,e}(\y_{b,c}\y_{d,f}-\y_{c,f}\y_{b,d})  
\\
-\l_{b}^2\y_{a,c}\y_{b,e}\y_{d,f}
)
\end{multline*}

\end{proof}

\begin{lem} The Riemannian curvature tensor applied to $\xi_{ij}$ is

\begin{multline*}
R_{ij,km}\xi_{ij}=
\frac{1}{8}\y_{j,m}\l_{j}^2(\l_{i}^2-3\l_{j}^2+3\l_{k}^2)\xi_{k,j}
+\frac{1}{8}\y_{i,k}\l_{i}^2(-3\l_{i}^2+\l_{j}^2+3\l_{m}^2)\xi_{i,m}
\\
-\frac{1}{4}\y_{i,m}\l_{i}^2(-\l_{j}^2 +\l_{i}^2-\l_{k}^2)\xi_{k,i}
+\frac{1}{8}\y_{i,m} (-\l_{i}^2+\l_{j}^2-\l_{k}^2)
(\l_{i}^2+\l_{j}^2 -\l_{k}^2 )\xi_{k,i}
\\
+ \frac{1}{4}\y_{j,k}\l_{j}^2(\l_{i}^2-\l_{j}^2+\l_{m}^2)\xi_{j,m} 
+\frac{1}{8}\y_{j,k}(-\l_{i}^2+\l_{j}^2+\l_{m}^2)
(-\l_{i}^2-\l_{j}^2+\l_{m}^2)\xi_{j,m}
\end{multline*}

\end{lem}

\begin{proof}
\begin{multline*}
R_{\xi_{ij} \xi_{km}}\xi_{ij}=R_{ij,km}\xi_{ij}=
\frac{1}{\sqrt{2}}\y_{j,k}\l_{j}^2\n_{im}\xi_{ij}
+\frac{1}{\sqrt{2}}\y_{j,m}\l_{j}^2\n_{(\xi_{ki}-\xi_{ik})}\xi_{ij}
\\
+\frac{1}{\sqrt{2}}\y_{i,k}\l_{i}^2\n_{(\xi_{mj}-\xi_{jm})}\xi_{ij}
-\frac{1}{\sqrt{2}}\y_{i,m}\l_{i}^2\n_{kj}\xi_{ij}
-\n_{ij}\n_{km}\xi_{ij}
=
\end{multline*}

\begin{multline*}
\frac{1}{4}\y_{j,k}\l_{i}^2\l_{j}^2\xi_{j,m}
-\frac{1}{4}\y_{j,k}\l_{j}^4\xi_{j,m}
+\frac{1}{4}\y_{j,k}\l_{j}^2\l_{m}^2\xi_{j,m}
\\ 
+\frac{1}{\sqrt{2}}\y_{j,m}\l_{j}^2\n_{(\xi_{ki}-\xi_{ik})}\xi_{ij}
+\frac{1}{\sqrt{2}}\y_{i,k}\l_{i}^2\n_{(\xi_{mj}-\xi_{jm})}\xi_{ij}
-\frac{1}{\sqrt{2}}\y_{i,m}\l_{i}^2\n_{kj}\xi_{ij}
\\
-
\frac{1}{2\sqrt{2}}(
+\l_{m}^2\y_{i,m}\n_{ij}\xi_{k,j}
+ \l_{m}^2\y_{m,j}\n_{ij}(\xi_{i,k}-\xi_{i,k})
+ \l_{k}^2\y_{i,k}\n_{ij}(\xi_{j,m}-\xi_{m,j}) 
\\ 
-\l_{k}^2\y_{k,j}\n_{ij}\xi_{i,m}
-\l_{j}^2\y_{i,k}\n_{ij}\xi_{j,m}
\\
-\l_{j}^2(\y_{i,m}\n_{ij}\xi_{k,j}-\y_{i,k}\n_{ij}\xi_{m,j})  
-\l_{i}^2(\y_{m,j}\n_{ij}\xi_{i,k}-\y_{k,j}\n_{ij}\xi_{i,m})  
+\l_{i}^2\y_{m,j}\n_{ij}\xi_{i,k}
\\
+\l_{k}^2\y_{m,j}\n_{ij}\xi_{i,k}
+\l_{m}^2(\y_{i,k}\n_{ij}\xi_{j,m}-\y_{k,j}\n_{ij}\xi_{i,m}) 
\\
+ \l_{k}^2(\y_{i,m}\n_{ij}\xi_{k,j}-\y_{m,j}\n_{ij}\xi_{i,k})  
-\l_{m}^2\y_{i,k}\n_{ij}\xi_{m,j})
.
\end{multline*}

\end{proof}

\begin{lem} The sectional curvature is

\begin{multline*}
K(\xi_{ij},\xi_{km})=(R_{ij,km}\xi_{ij},\xi_{km})=
\\
\frac{1}{8}\y_{j,m}\l_{j}^2(\l_{i}^2-3\l_{j}^2+3\l_{k}^2)
+\frac{1}{8}\y_{i,k}\l_{i}^2(-3\l_{i}^2+\l_{j}^2+3\l_{m}^2)
\\
+\frac{1}{4}\y_{i,m}\l_{i}^2(\l_{j}^2-\l_{i}^2+\l_{k}^2)
+\frac{1}{8}\y_{i,m} (-\l_{i}^2+\l_{j}^2-\l_{k}^2)
(\l_{i}^2+\l_{j}^2 -\l_{k}^2 )
\\
+ \frac{1}{4}\y_{j,k}\l_{j}^2(\l_{i}^2-\l_{j}^2+\l_{m}^2) 
+\frac{1}{8}\y_{j,k}(-\l_{i}^2+\l_{j}^2+\l_{m}^2)
(-\l_{i}^2-\l_{j}^2+\l_{m}^2)
\end{multline*}

\end{lem}

The third part of the following theorem says that the principal Ricci curvatures for $\mathfrak{s}\mathfrak{o}_{HS}$ can tend to either $\infty$ or $-\infty$ as the dimension $N \to \infty$ depending on the choice of the scaling $\{\l_i\}_{i=1}^{\infty}$. 

\begin{thm}\ll{T:O} Let $N>\max{i,j}$. Then 

\begin{enumerate} 

\item the truncated Ricci curvature is

\[
R_{ij}^N=
\frac{1}{8}(-4\l_{i}^2\l_{j}^2-5(\l_{i}^2-\l_{j}^2)^2)N
+\frac{3}{8}(\l_{i}^2+\l_{j}^2)\sum_{m=1}^{m=N}\l_{m}^2
+\frac{1}{4}\sum_{m=1}^{m=N}\l_{m}^4.
\]

\item For the Hilbert-Schmidt inner product the truncated Ricci curvature is

\[
R_{ij}^N=\frac{N}{2}
.
\]

\item
The truncated self-adjoint Ricci curvature is diagonal in the basis $\{\xi_{km}\}$. Let $\{a_{km}\}$ be its principal Ricci curvatures, that is, $\hat R^N(\xi_{km})=a_{km}\xi_{km}$. Then if $\{\l_i\}_{i=1}^{\infty} \in \ell^2$, the principal Ricci curvatures have the following asymptotics as $N \to \infty$

\[
a_{km}=
A_{km}^N
+\frac{N}{8}
(2\l_{k}^2-3\l_{m}^2)(\l_{m}^2-\l_{k}^2), 
\]
where $A_{km}^N\to A<\infty$ as $N \to \infty$.

\end{enumerate}
\end{thm}

\begin{proof}

(1)
\begin{multline*}
R_{ij}^N=\sum_{k,m=1}^{k,m=N}(R_{ij,km}\xi_{ij}, \xi_{km})
=\\
\sum_{k=1}^{k=N}
\frac{1}{8}\l_{j}^2(\l_{i}^2-3\l_{j}^2+3\l_{k}^2)
+\frac{1}{8}\y_{i,k}\l_{i}^2(-3\l_{i}^2N+\l_{j}^2N
+3\sum_{m=1}^{m=N}\l_{m}^2)
\\
+\frac{1}{4}\l_{i}^2(\l_{j}^2-\l_{i}^2+\l_{k}^2)
+\frac{1}{8}(-\l_{i}^2+\l_{j}^2-\l_{k}^2)(\l_{i}^2+\l_{j}^2 -\l_{k}^2 )
\\
+ \frac{1}{4}\y_{j,k}\l_{j}^2(\l_{i}^2N-\l_{j}^2N+\sum_{m=1}^{m=N}\l_{m}^2)
\\ 
+\frac{1}{8}\y_{j,k}((\l_{i}^2-\l_{j}^2)(\l_{i}^2+\l_{j}^2)N-2\l_{i}^2\sum_{m=1}^{m=N}\l_{m}^2+\sum_{m=1}^{m=N}\l_{m}^4)
=
\\
+\frac{1}{8}(6\l_{i}^2\l_{j}^2-5\l_{i}^4-5\l_{j}^4)N
+\frac{3}{8}(\l_{i}^2+\l_{j}^2)\sum_{m=1}^{m=N}\l_{m}^2
+\frac{1}{4}\sum_{m=1}^{m=N}\l_{m}^4.
\end{multline*}

(2) follows from (1).

(3)
Let $k<m<N$, then  the truncated self-adjoint Ricci curvature is

\begin{multline*}
\hat R^N(\xi_{km})=\sum_{i<j\leqslant N}R_{ij,km}\xi_{ij}
=\sum_{j=1}^{N}\sum_{i=1}^{j-1}R_{ij,km}\xi_{ij}
=
\\
\sum_{j=1}^{k}\sum_{i=1}^{j-1}R_{ij,km}\xi_{ij}
+\sum_{j=k+1}^{m}\sum_{i=1}^{j-1}R_{ij,km}\xi_{ij}
+\sum_{j=m+1}^{N}\sum_{i=1}^{j-1}R_{ij,km}\xi_{ij}.
\end{multline*}
Thus

\begin{multline*}
\hat R^N(\xi_{km})=
\frac{1}{8}\sum_{j=1}^{k}\sum_{i=1}^{j-1} \y_{j,m}\l_{j}^2(\l_{i}^2-3\l_{j}^2+3\l_{k}^2)\xi_{k,j}
\\
+2\y_{j,k}\l_{j}^2(\l_{i}^2-\l_{j}^2+\l_{m}^2)\xi_{j,m} 
+\y_{j,k}(-\l_{i}^2+\l_{j}^2+\l_{m}^2)
(-\l_{i}^2-\l_{j}^2+\l_{m}^2)\xi_{j,m}
\\
+\sum_{j=k+1}^{m}\sum_{i=1}^{j-1}R_{ij,km}\xi_{ij}
+\sum_{j=m+1}^{N}\sum_{i=1}^{j-1}R_{ij,km}\xi_{ij}
=
\end{multline*}

\begin{multline*}
\frac{1}{8}
\left[
\sum_{i=1}^{k-1}
2\l_{k}^2(\l_{i}^2-\l_{k}^2+\l_{m}^2)
+(-\l_{i}^2+\l_{k}^2+\l_{m}^2)(-\l_{i}^2-\l_{k}^2+\l_{m}^2)
\right.
\\
+\l_{m}^2(\l_{k}^2-3\l_{m}^2+3\l_{k}^2)
+\sum_{j=k+1}^{m}\l_{k}^2(-3\l_{k}^2+\l_{j}^2+3\l_{m}^2)
\\
+\sum_{i=1, i \not= k}^{m-1}
\l_{m}^2(\l_{i}^2-3\l_{m}^2+3\l_{k}^2)
+\sum_{j=m+1}^{N}
\l_{k}^2(-3\l_{k}^2+\l_{j}^2+3\l_{m}^2)
\\
\left.
+\sum_{j=m+1}^{N}
-2\l_{m}^2(-\l_{j}^2 +\l_{m}^2-\l_{k}^2)
+(-\l_{m}^2+\l_{j}^2-\l_{k}^2)
(\l_{m}^2+\l_{j}^2 -\l_{k}^2 )
\right]\xi_{k,m}
=
\end{multline*}

\begin{multline*}
\frac{1}{8}
\left[
\sum_{l=1}^{k-1}
2\l_{k}^2(\l_{l}^2-\l_{k}^2+\l_{m}^2) 
+(\l_{l}^2-\l_{k}^2-\l_{m}^2)
(\l_{l}^2+\l_{k}^2-\l_{m}^2)
\right.
\\
+\sum_{l=k+1}^{N}
\l_{k}^2(\l_{l}^2-3\l_{k}^2+3\l_{m}^2)
+\sum_{l=1}^{m-1}
\l_{m}^2(\l_{l}^2+3\l_{k}^2-3\l_{m}^2)
\\
\left.
+\sum_{l=m+1}^{N}
-2\l_{m}^2(-\l_{l}^2 -\l_{k}^2+\l_{m}^2)
+(\l_{l}^2-\l_{k}^2-\l_{m}^2)
(\l_{l}^2 -\l_{k}^2 +\l_{m}^2)
\right]
\xi_{k,m}
=
\end{multline*}

\begin{multline*}
\frac{1}{8}
[
(k-1)(3\l_{k}^2+\l_{m}^2)
-3(m-1)\l_{m}^2
(-\l_{k}^2+\l_{m}^2)
\\
+
(N-k)
(\l_{k}^2(-3\l_{k}^2+3\l_{m}^2))
+(N-m)
(-\l_{k}^2-3\l_{m}^2)
(-\l_{k}^2 +\l_{m}^2)
\\
+
\sum_{l=1}^{k-1}
(2\l_{k}^2\l_{l}^2
+\l_{l}^4-2\l_{m}^2\l_{l}^2)
+\sum_{l=k+1}^{N}
\l_{k}^2\l_{l}^2
\\
+\sum_{l=1}^{m-1}
\l_{m}^2\l_{l}^2
+\sum_{l=m+1}^{N}
\l_{l}^2(2\l_{m}^2 +\l_{l}^2-2\l_{k}^2)
]
\xi_{k,m}.
\end{multline*}

\end{proof}

\subsection{Upper triangular Hilbert-Schmidt algebra} As before let $\{\l_i\}_{i=1}^{\infty}$ be a bounded sequence of strictly positive numbers. 
In this section we consider the infinite dimensional Lie algebra $\G$ generated by the orthonormal basis $\xi_{ij}=\l_i\l_je_{ij}$, $i<j$. In this case $\G$ is a Lie subalgebra of $\mathfrak{h}_{HS}$. Recall that if the sequence $\{\l_i\}_{i=1}^{\infty}$ is bounded, then by Theorem \ref{T:split}, Continuity Assumption \ref{E:norm} is satisfied for the corresponding  norm on $\G$.

The Lie bracket is
\[
[\xi_{ij},\xi_{km}]=\y_{j,k}\l_{j}^2\xi_{im}-
\y_{i,m}\l_{i}^2\xi_{kj}, \ i<j, \ , k<m.
\]
Denote $\n_{ab}=\n_{\xi_{ab}}$, then as for the general algebra with the convention that $\xi_{ij}=0$ if $i\geqslant j$

\[
\n_{ab}\xi_{cd}
=
\frac 12\left(
\y_{b,c}\l_{b}^2\xi_{ad}
-\y_{a,d}\l_{a}^2\xi_{cb}
-\y_{a,c}\l_{d}^2\xi_{db}
+\y_{b,d}\l_{c}^2\xi_{ac}
+\y_{b,d}\l_{a}^2\xi_{ca}
-\y_{a,c}\l_{b}^2\xi_{bd}
\right).
\]

\begin{lem}

The Riemannian curvature tensor applied to $\xi_{ij}$ is

\begin{multline*}
R_{\xi_{ij} \xi_{km}}\xi_{ij}
=
\frac 14(
2\y_{i,k}\y_{j,m}\l_{i}^4\xi_{ij}
+2\y_{i,k}\y_{j,m}\l_{j}^4\xi_{ij}
-3\y_{i,m}\l_{i}^4\xi_{ki}
-3\y_{j,k}\l_{j}^4\xi_{jm}
\\
+\y_{i,k}\l_{j}^4\xi_{im}
+\y_{i,k}\l_{m}^4\xi_{im}
+\y_{j,m}\l_{i}^4\xi_{kj}
+\y_{j,m}\l_{k}^4\xi_{kj}
)
\end{multline*}

The sectional curvature is
\begin{multline*}
(R_{ij,km}\xi_{ij}, \xi_{km})
=\\
\frac 14(
2\y_{i,k}\y_{j,m}\l_{i}^4   
+2\y_{i,k}\y_{j,m}\l_{j}^4   
-3\y_{i,m}\l_{i}^4
-3\y_{j,k}\l_{j}^4
+\y_{i,k}\l_{j}^4
+\y_{i,k}\l_{m}^4
+\y_{j,m}\l_{i}^4
+\y_{j,m}\l_{k}^4
)
.
\end{multline*}

\end{lem}

\begin{proof}
The Riemannian curvature tensor applied to $\xi_{ij}$ is

\begin{multline*}
R_{\xi_{ij} \xi_{km}}\xi_{ij}=R_{ij,km}\xi_{ij}
=\\
\y_{j,k}\l_{j}^2\n_{im}\xi_{ij}-\y_{i,m}\l_{i}^2\n_{kj}\xi_{ij}
-\n_{ij}\n_{km}\xi_{ij}+
\n_{km}\n_{ij}\xi_{ij}
=\\
\frac 12(
-\y_{j,k}\l_{j}^4\xi_{jm}
-\y_{i,m}\l_{i}^4\xi_{ki}
)
+\frac 14(
+2\y_{i,k}\y_{j,m}\l_{i}^4\xi_{ij}
+2\y_{i,k}\y_{j,m}\l_{j}^4\xi_{ij}
\\
-\y_{i,m}\l_{i}^4\xi_{ki}
-\y_{j,k}\l_{j}^4\xi_{jm}
+\y_{i,k}\l_{j}^4\xi_{im}
+\y_{i,k}\l_{m}^4\xi_{im}
+\y_{j,m}\l_{i}^4\xi_{kj}
+\y_{j,m}\l_{k}^4\xi_{kj}
)
=\\
\frac 14(
2\y_{i,k}\y_{j,m}\l_{i}^4\xi_{ij}
+2\y_{i,k}\y_{j,m}\l_{j}^4\xi_{ij}
-3\y_{i,m}\l_{i}^4\xi_{ki}
-3\y_{j,k}\l_{j}^4\xi_{jm}
\\
+\y_{i,k}\l_{j}^4\xi_{im}
+\y_{i,k}\l_{m}^4\xi_{im}
+\y_{j,m}\l_{i}^4\xi_{kj}
+\y_{j,m}\l_{k}^4\xi_{kj}
).
\end{multline*}

\end{proof}

The third part of the following theorem says that the principal Ricci curvatures for $\mathfrak{h}_{HS}$ tend to $-\infty$ as the dimension $N \to \infty$. This can interpreted as the Ricci curvature being $-\infty$. Note that the condition on the scaling $\{\l_i\}_{i=1}^{\infty} \in \ell^2$ corresponds to the condition we assumed in \cite{Gor00-1}, \cite{Gor00-2}, \cite{Gor02}, \cite{Gor03-2}. We needed this condition to construct a heat kernel measure living in $HS+I$. This means that such a measure exists even though the Ricci curvature is $-\infty$.

\begin{thm}\ll{T:T} 
Let $N>\max{i,j}$ then

\begin{enumerate}
\item the truncated Ricci curvature is

\begin{multline*}
R_{ij}^N=\sum_{k,m=1}^{k,m=N}(R_{ij,km}\xi_{ij}, \xi_{km})
=\\
\frac 14(
(4-3i+j)\l_{i}^4
+(2+3j-i-2N)\l_{j}^4
+\sum_{l=i+1}^{N}\l_{l}^4
+\sum_{l=1}^{j-1}\l_{l}^4
)
.
\end{multline*}
For the Hilbert-Schmidt inner product the truncated Ricci curvature 

\[
R_{ij}^N=
\frac 14(5-5i+5j-N)
.
\]

\item
The truncated self-adjoint Ricci curvature is diagonal in the basis $\{\xi_{km}\}$. Let $\{a_{km}\}$ be its principal Ricci curvatures, that is, $\hat R^N(\xi_{km})=a_{km}\xi_{km}$. Then if $\{\l_i\}_{i=1}^{\infty} \in \ell^2$, the principal Ricci curvatures have the following asymptotics as $N \to \infty$

\[
a_{km}=
\left(
B_{km}^N
-\frac N2 \l_{m}^4
\right)
\xi_{k,m}, 
\]
where $B_{km}^N\to B_{km}<\infty$ as $N \to \infty$. 
\end{enumerate}
\end{thm}

\begin{proof}
(1) follows directly from the previous lemmas. 
(2) Let $k<m<N$, then the truncated adjoint Ricci curvature is

\begin{multline*}
\hat R^N(\xi_{km})=\sum_{i<j\leqslant N}R_{ij,km}\xi_{ij}
=\sum_{j=1}^{N}\sum_{i=1}^{j-1}R_{ij,km}\xi_{ij}
=\\
\sum_{j=1}^{k}\sum_{i=1}^{j-1}R_{ij,km}\xi_{ij}
+\sum_{j=k+1}^{m}\sum_{i=1}^{j-1}R_{ij,km}\xi_{ij}
+\sum_{j=m+1}^{N}\sum_{i=1}^{j-1}R_{ij,km}\xi_{ij}.
\end{multline*}

\begin{multline*}
\hat R^N(\xi_{km})=\sum_{i<j\leqslant N}R_{ij,km}\xi_{ij}
=\sum_{j=1}^{N}\sum_{i=1}^{j-1}R_{ij,km}\xi_{ij}
=\\
\sum_{j=1}^{k}\sum_{i=1}^{j-1}R_{ij,km}\xi_{ij}
+\sum_{j=k+1}^{m}\sum_{i=1}^{j-1}R_{ij,km}\xi_{ij}
+\sum_{j=m+1}^{N}\sum_{i=1}^{j-1}R_{ij,km}\xi_{ij}
=\\
\frac {(-3k+m+4)}4
\l_{k}^4\xi_{km}
+\frac {(2-k-2N+3m)}4
\l_{m}^4\xi_{km}
+\frac 14\sum_{l=k+1}^{N}\l_{l}^4\xi_{km}
+\frac 14\sum_{l=1}^{m-1} \l_{l}^4\xi_{km}
=\\
B_{km}^N\xi_{km}
-\frac N2 \l_{m}^4\xi_{km}.
\end{multline*}

\end{proof}

\begin{cor}
Suppose $\{\l_i\}_{i=1}^{\infty} \in \ell^2$. Then for large $N$ 

\[
\hat R^N(\xi_{km})=b_{km }^N \xi_{km}, 
\]
where $b_{km }^N\to -\infty$ as $N\to \infty$. This can be described as the Ricci curvature being negative infinity for any $x \in \g$.
\end{cor}


\begin{thebibliography}{99}


\bibitem{ARS85}
M. Adams, T. Ratiu, R.  Schmid,
\emph{The Lie group structure of diffeomorphism groups and invertible Fourier integral operators, with applications},
Infinite-dimensional groups with applications (Berkeley, Calif., 1984), 1--69,
Math. Sci. Res. Inst. Publ., \textbf{4}, Springer, New York, 1985.

\bibitem{Ci89}
R. Cianci, 
\emph{Infinitely generated supermanifolds, super Lie groups and homogeneous supermanifolds},
Symmetries in science, III (Vorarlberg, 1988), 147--162,
Plenum, New York, 1989.

\bibitem{Cz89}
J. Czyz, 
\emph{On Lie supergroups and superbundles defined via the Baker-Campbell-Hausdorff formula},
J. Geom. Phys. \textbf{6}, 1989, no. 4, 595--626.

\bibitem{Har72}
Pierre de la Harpe, 
\emph{Classical Banach-Lie algebras and Banach-Lie groups of operators in Hilbert space},
Lecture Notes in Mathematics, \textbf{285}, 
Springer-Verlag, Berlin-New York, 1972.

\bibitem{Dr92}  
B. Driver,
\emph{A Cameron-Martin type quasi-invariance theorem for Brownian motion on a compact Riemannian manifold},
J. Funct. Anal. \textbf{110}, 1992, no. 2, 272--376.

\bibitem{Dr94}  
B. Driver, \emph{A Cameron-Martin type quasi-invariance
theorem for pinned Brownian motion on a compact Riemannian manifold}, Trans.
Amer. Math. Soc., \textbf{342}, 1994, no. 1, 375-395.

\bibitem{Driver95}
B. Driver, 
\emph{Towards calculus and geometry on path spaces},
Stochastic analysis (Ithaca, NY, 1993), Proc. Sympos. Pure Math.,
\textbf{57}, 1995, pp. 405--422, Amer. Math. Soc., Providence, RI.

\bibitem{DriverLohr}
B. Driver, T. Lohrenz,
\emph{Logarithmic Sobolev inequalities for pinned loop groups}, J. Funct. Anal.,
 \textbf{140}, 1996, pp. 381--448.

\bibitem{DriverHu}
B. Driver, Y. Hu, 
\emph{On heat kernel logarithmic Sobolev inequalities},
 booktitle={Stochastic analysis and applications (Powys, 1995)}, pp. 189--200, World Sci. Publishing, River Edge, NJ, 1996.

\bibitem{Dr97}  
B. Driver, 
\emph{Integration by parts and quasi-invariance
for heat kernel measures on loop groups}, J. of Funct. Anal., \textbf{149},
1997, no. 2, 470-547.

\bibitem{Dr97-2}  
B. Driver,
\emph{Integration by parts for heat kernel measures revisited},
J. Math. Pures Appl. (9) \textbf{76}, 1997, no. 8, 703--737.

\bibitem{Driver98}
B. Driver, 
\emph{A correction to the paper: ``Integration by parts and quasi-invariance for heat kernel measures on loop groups''}, J. Funct. Anal.,
\textbf{155}, 1998, pp. 297--301.

\bibitem{Driver03}
B. Driver, 
\emph{Analysis of Wiener measure on path and loop groups},
 Finite and infinite-dimensional analysis in honor of Leonard Gross (New Orleans, LA, 2001),
Contemp. Math., \textbf{317}, pp. 57--85, Amer. Math. Soc., Providence, RI, 2003.

\bibitem{Dyn47}
E. B. Dynkin, 
\emph{Calculation of the coefficients in the Campbell-Hausdorff formula},
Doklady Akad. Nauk SSSR, \textbf{57}, 1947, 323--326.

\bibitem{Dyn49}
E. B. Dynkin, 
\emph{On the representation by means of commutators of the series ${\rm log} (e^x e^y)$ for noncommutative $x$ and $y$}, Mat. Sbornik N.S., \textbf{25(67)}, 1949, 155--162.

\bibitem{Fang04}
S. Fang, 
\emph{Metrics $\H_s$ and behaviours as $s\downarrow 1/2$ on loop groups}, J. Funct. Anal., \textbf{213}, 2004, 440--465.

\bibitem{Fang99}
S. Fang, 
\emph{Ricci tensors on some infinite-dimensional Lie algebras}, 
J. Funct. Anal., \textbf{161}, 1999, pp. 132--151.

\bibitem{Freed88}
D. Freed, 
\emph{The geometry of loop groups}, 
 J. Differential Geom., \textbf{28}, 1988, pp. 223--276.

\bibitem{FukushimaBook}
M. Fukushima, 
\emph{Dirichlet forms and Markov processes}, 
1980, North-Holland Mathematical Library, \textbf{23}.

\bibitem{Gl02-1}
H. Gl\"{o}ckner, 
\emph{Algebras whose groups of units are Lie groups},
Studia Math. \textbf{153}, 2002, no. 2, 147--177.

\bibitem{Gl02-2}
H.Gl\"{o}ckner,
\emph{Infinite-dimensional Lie groups without completeness restrictions},
Geometry and analysis on finite- and infinite-dimensional Lie groups, 43--59,
Banach Center Publ., \textbf{55},
Polish Acad. Sci., Warsaw, 2002.

\bibitem{Gl02-3}
H. Gl\"{o}ckner, 
\emph{Lie group structures on quotient groups and universal complexifications for infinite-dimensional Lie groups},
J. Funct. Anal. \textbf{194}, 2002, no. 2, 347--409.

\bibitem{Gl03}
H. Gl\"{o}ckner,
\emph{Direct limit Lie groups and manifolds}, 
J. Math. Kyoto Univ. \textbf{43}, 2003, no. 1, 2--26.

\bibitem{Gor00-1} 
M. Gordina,
\emph{Holomorphic functions and the heat kernel
measure on an infinite-dimensional complex orthogonal group}, Potential
Analysis Volume, \textbf{12}, 2000, pp. 325-357.

\bibitem{Gor00-2}  
M. Gordina,
\emph{Heat kernel analysis and Cameron-Martin
subgroup for infinite-dimensional groups}, J. Func. Anal., \textbf{171},
2000, pp. 192-232.

\bibitem{Gor02}  
M. Gordina,
\emph{Taylor map on groups associated with a
$\text{II}_1$-factor}, 2002, Infin. Dimens. Anal. Quantum Probab. Relat. Top.,
\textbf{5}, 93--111.

\bibitem{Gor03-2} 
M. Gordina,
\emph{
Stochastic differential equations on noncommutative $L^2$
}, 2003,
Contemp. Math., Amer. Math. Soc.,  Providence, RI
\textbf{317}, 87--98.

\bibitem{HsuBook}
E. Hsu,
\emph{Stochastic analysis on manifolds},
Graduate Studies in Mathematics,
2002,
American Mathematical Society,\textbf{38}.

\bibitem{HsuIAS}
E. Hsu,
\emph{Analysis on path and loop spaces},
Probability theory and applications (Princeton, NJ, 1996),
277--347, IAS/Park City Math. Ser. \textbf{6}, Amer. Math. Soc., 1999.

\bibitem{Hsu99}
E. Hsu,
\emph{Estimates of derivatives of the heat kernel on a compact Riemannian manifold},
Proc. Amer. Math. Soc. \textbf{127}, 1999, no. 12, 3739--3744.

\bibitem{Hsu95}
E. Hsu,
\emph{Quasi-invariance of the Wiener measure on the path space over a compact Riemannian manifold},
J. Funct. Anal. \textbf{134}, 1995, no. 2, 417--450.

\bibitem{Hsu93}
E. Hsu,
\emph{Flows and quasi-invariance of the Wiener measure on path spaces},
Stochastic analysis (Ithaca, NY, 1993), 265--279, Proc. Sympos. Pure Math., \textbf{57}, Amer. Math. Soc., 1995.

\bibitem{KM97}
A. Kriegl, P. Michor,
\emph{Regular infinite-dimensional Lie groups},
J. Lie Theory \textbf{7}, 1997, no. 1, 61--99.

\bibitem{LT65}
M. Lazard, J. Tits, 
\emph{Domaines d'injectivit\'{e} de l'application exponentielle} (French), 
Topology \textbf{4}, 1965/1966, 315--322.

\bibitem{L95}
L. Lempert, 
\emph{The problem of complexifying a Lie group},
Multidimensional complex analysis and partial differential equations (S\~{a}o Carlos, 1995), 169--176,
Contemp. Math., \textbf{205},
Amer. Math. Soc., Providence, RI, 1997.

\bibitem{Le92}
J. Leslie,
\emph{Some integrable subalgebras of the Lie algebras of infinite-dimensional Lie groups},
Trans. Amer. Math. Soc. \textbf{333}, 1992, no. 1, 423--443.

\bibitem{M01}
P. Malliavin, 
\emph{Probability and geometry},
Taniguchi Conference on Mathematics Nara '98, 179--209,
Adv. Stud. Pure Math., \textbf{31},
Math. Soc. Japan, Tokyo, 2001. 


\bibitem{MR95}
J. Marion, T. Robart,
\emph{Regular Fr\'{e}chet-Lie groups of invertible elements in some inverse limits of unital involutive Banach algebras}, Georgian Math. J. \textbf{2}, 1995, no. 4, 425--444.

\bibitem{Milnor76} 
J. Milnor,
\emph{Curvatures of left invariant metrics on {L}ie groups},
Advances in Math.,
\textbf{21}, 1976, pp. 293--329.

\bibitem{Mi84}
J. Milnor, 
\emph{Remarks on infinite-dimensional Lie groups},
Relativity, groups and topology, \textbf{II} (Les Houches, 1983), 1007--1057,
North-Holland, Amsterdam, 1984.

\bibitem{NRW01}
L. Natarajan, E. Rodr\'{i}guez-Carrington, J. Wolf, 
\emph{The Bott-Borel-Weil theorem for direct limit groups}, 
Trans. Amer. Math. Soc. \textbf{353}, 2001, no. 11, 4583--4622.

\bibitem{NRW91}
L. Natarajan, E. Rodr\'{i}guez-Carrington, J. Wolf, 
\emph{Differentiable structure for direct limit groups},
Lett. Math. Phys. \textbf{23}, 1991, no. 2, 99--109.

\bibitem{NRW94}
L. Natarajan, E. Rodr\'{i}guez-Carrington, J. Wolf, 
\emph{New classes of infinite-dimensional Lie groups},
Algebraic groups and their generalizations: quantum and infinite-dimensional methods (University Park, PA, 1991), 377--392,
Proc. Sympos. Pure Math., \textbf{56}, Part 2,
Amer. Math. Soc., Providence, RI, 1994. 

\bibitem{NRW93}
L. Natarajan, E. Rodr\'{i}guez-Carrington, J. Wolf, 
\emph{Locally convex Lie groups},
Nova J. Algebra Geom. \textbf{2}, 1993, no. 1, 59--87.

\bibitem{Neeb98}
Karl-Hermann Neeb, 
\emph{Holomorphic highest weight representations of infinite-dimensional complex classical groups},
J. Reine Angew. Math. \textbf{497}, 1998, 171--222.

\bibitem{Ne96}
Yu. A. Neretin, 
\emph{Categories of symmetries and infinite-dimensional groups},
Translated from the Russian by G. G. Gould. London Mathematical Society Monographs. New Series, \textbf{16},  Oxford Science Publications, 
The Clarendon Press, Oxford University Press, New York, 1996. 

\bibitem{OMYK83}
H. Omori, Y. Maeda, A. Yoshioka, O. Kobayashi, 
\emph{On regular Fr\'{e}chet-Lie groups. V. Several basic properties},
Tokyo J. Math. \textbf{6}, 1983, no. 1, 39--64.

\bibitem{SegalBook} 
A. Pressley, G. Segal,
\emph{Loop groups}, 
Oxford Mathematical Monographs, Oxford University, 1986.

\bibitem{R96}
T. Robart, 
\emph{Groupes de Lie de dimension infinie. Second et troisi\`{e}me th\`{e}or\'{e}mes de Lie. I. Groupes de premi\`{e}re esp\`{e}ce}, (French. English, French summary) [Infinite-dimensional Lie groups. Second and third theorems. I. Groups of the first kind]
C. R. Acad. Sci. Paris S\'{e}r. I Math. \textbf{322}, 1996, no. 11, 1071--1074.

\bibitem{R97}
T. Robart,
\emph{Sur l'int\'{e}grabilit\'{e} des sous-alg\`{e}bres de Lie en dimension infinie}, (French. French summary) [On the integrability of infinite-dimensional Lie subalgebras]
Canad. J. Math. \textbf{49}, 1997, no. 4, 820--839.

\bibitem{R02}
T. Robart,
\emph{Around the exponential mapping},
Infinite dimensional Lie groups in geometry and representation theory (Washington, DC, 2000),
World Sci. Publishing, River Edge, NJ, 2002, 11--30.

\bibitem{Shig97}
I. Shigekawa, 
\emph{Differential calculus on a based loop group}, in "New Trends in Stochastic analysis" (K. D. Elworthy, S. Kusuoka, I. Shigekawa, Eds.), pp. 375--398, 1997.

\bibitem{SperaWurz00}
M. Spera,
T. Wurzbacher, 
\emph{Differential geometry of Grassmannian embeddings of based loop groups},
Differential Geom. Appl., \textbf{13}, 2000, pp. 43--75.

\bibitem{ShigTanig} 
I. Shigekawa, S. Taniguchi, 
\emph{A K\"ahler metric on a based loop group and a covariant differentiation}, 
in It\^o's stochastic calculus and probability theory, Springer, Tokyo,
1996, pp. 327-346.

\bibitem{StroockBook} 
D. Stroock, 
\emph{An introduction to the analysis of paths on a Riemannian manifold}, 
Mathematical Surveys and Monographs, American Mathematical Society, 
Providence, RI, 2000.

\bibitem{Tan96}
S. Taniguchi,
\emph{On Ricci curvatures of hypersurfaces in abstract Wiener spaces},
J. Funct. Anal., \textbf{136}, 1996, 226--244.

\bibitem{V72}
F.-H. Vasilescu,
\emph{Normed Lie algebras},
Canad. J. Math. \textbf{24}, 1972, 580--591.

\bibitem{W98}
W. Wojty\'nski, 
\emph{Quasinilpotent Banach-Lie algebras are Baker-Campbell-Hausdorff}, J. Funct. Anal. \textbf{153}, 1998, no. 2, 405--413.

\end{thebibliography}
\end{document}